\documentclass[12pt]{article}
\author{S\'everine Fiedler-Le Touz\'e}
\usepackage{amsfonts, amssymb, epsfig}
\newtheorem{theorem}{Theorem}
\newtheorem{conjecture}{Conjecture}
\newtheorem{definition}{Definition}
\newtheorem{lemma}{Lemma}
\newtheorem{proposition}{Proposition}
\title{M-curves of degree 9 with three nests}
\begin{document}
\maketitle
\begin{abstract}
The first part of Hilbert's sixteenth problem deals with the classification of the  
isotopy types realizable by real plane algebraic curves of a given degree $m$.
For $m = 9$ the classification of the $M$-curves is still wide open.
Let $C_9$ be an $M$-curve of degree $9$ and $O$ be a non-empty oval of $C_9$. If $O$ contains in its interior
$\alpha$ ovals that are all empty, we say that $O$ together with these $\alpha$ ovals forms a nest.
The present paper deals with the $M$-curves with three nests. Let $\alpha_i, i=1, 2, 3$ be the numbers of
empty ovals in each nest. We prove that at least one of the $\alpha_i$ is odd. 
This excludes $41$ new isotopy types. 
This is a step towards a conjecture of A. Korchagin, claiming that at least two of the $\alpha_i$ should be odd. 
\end{abstract}

\section{Introduction}

\subsection{Real and complex schemes}
The first part of Hilbert's sixteenth problem deals with the classification
of the isotopy types realizable by real plane algebraic curves of given degree.
Let $A$ be a real algebraic non-singular plane curve of degree $m$.
Its complex part $\mathbb{C}A \subset \mathbb{C}P^2$ is a Riemannian surface
of genus $g=(m-1)(m-2)/2$; its real part $\mathbb{R}A \subset \mathbb{R}P^2$
is a collection of $L \leq g+1$ circles embedded in $\mathbb{R}P^2$.
If $L=g+1$, we say that A is an $M$-curve.A circle embedded in
$\mathbb{R}P^2$ is called {\em oval\/} or {\em pseudo-line\/} depending
on whether it realizes the class 0 or 1 of $H_1(\mathbb{R}P^2)$.If $m$ is even,
the $L$ components of $\mathbb{R}A$ are ovals; if $m$ is odd,    
$\mathbb{R}A$ contains exactly one pseudo-line,which will be denoted by
$\mathcal{J}$. An oval separates $\mathbb{R}P^2$  into a M\"obius band and
a disc. The latter is called the {\em interior\/} of the oval.
An oval of $\mathbb{R}A$ is {\em empty\/} if its interior contains no other
oval. One calls {\em exterior oval\/} an oval that is surrounded by no other
oval. 
Two ovals form an {\em injective pair\/} if one of them lies in the
interior of the other one.
Let us call the isotopy type of  $\mathbb{R}A \subset \mathbb{R}P^2$ the
{\em real scheme\/} of $A$; it will be described with the following
notation due to Viro. The symbol $\langle \mathcal{J} \rangle$
stands for a curve consisting of one single pseudo-line; 
$\langle n \rangle$ stands for a curve consisting of $n$ empty ovals.
If $X$ is the symbol for a curve without pseudo-line, $1 \langle X \rangle$
is the curve obtained by adding a new oval, containing all of the
others in its interior. Finally, a curve which is the union of two disjoint 
curves $\langle A \rangle$ and $\langle B \rangle$, having the property that
none of the ovals of one curve is contained in an oval of the other curve,
is denoted by $\langle A \amalg B \rangle$. The classification of the real 
schemes which are realizable by $M$-curves is complete up to degree 7, see for example \cite{vi1}, \cite{vi2}, \cite{wi}.
For $m \geq 8$, one restricts the study to the case of the $M$-curves. The classification is almost complete
for $m = 8$, and still wide open for $m = 9$.
A systematic study of the case $m = 9$ has been done, the main 
contribution being due to A. Korchagin. See e.g. \cite{ko6},  \cite{ko3}, 
\cite{ko4}, \cite{ko5}, \cite{or3}, \cite{or7}
for the constructions, and \cite{ko1}, \cite{ko2}, \cite {ko6}, \cite{fi}, \cite{flt1}, \cite{flt2}, \cite{flt3}
\cite{or2}, \cite{or4}, \cite{or5}, \cite{or7} 
for the restrictions.

Let us briefly recall some facts about complex orientations, see also \cite{ro}, \cite{fi}.
The complex conjugation $conj$ of $\mathbb{C}P^2$ acts on $\mathbb{C}A$
with $\mathbb{R}A$ as fixed points sets. Thus, $\mathbb{C}A \setminus
\mathbb{R}A$ is connected, or splits in 2 homeomorphic halves which are
exchanged by $conj$. In the latter case, we say that $A$ is dividing.
Let us now consider a dividing curve $A$ of degree $m$, and
assume that $\mathbb{C}A$ is oriented canonically.
We choose a half $\mathbb{C}A_+$ of $\mathbb{C}A \setminus \mathbb{R}A$.
The orientation of $\mathbb{C}A_+$ induces an orientation on its
boundary $\mathbb{R}A$. This orientation, which is defined up to complete
reversion, is called {\em complex orientation\/} of $A$.
One can provide all the injective pairs of $\mathbb{R}A$ with a sign as
follows: such a pair is {\em positive\/} if and only if the orientations
of its 2 ovals induce an orientation of the annulus that they bound in
$\mathbb{R}P^2$. Let $\Pi_+$ and $\Pi_-$ be the numbers of positive
and negative injective pairs of $A$. If $A$ has odd degree, each
oval of $\mathbb{R}A$ can be provided with a sign: given an oval $O$
of $\mathbb{R}A$, consider the M\"obius band $\mathcal{M}$
obtained by cutting away the interior of $O$ from $\mathbb{R}P^2$.
The classes $[O]$ and $[2\mathcal{J}]$ of $H_1(\mathcal{M})$ either
coincide or are opposite. In the first case, we say that $O$ is
{\em negative\/}; otherwise $O$ is {\em positive\/}. Let $\Lambda_+$ and $\Lambda_-$
be respectively the numbers of positive and negative ovals of
$\mathbb{R}A$. 
The {\em complex scheme\/} of $A$ is obtained by
enriching the real scheme with the complex orientation:
let e.g. $A$ have real scheme $\langle \mathcal{J} \amalg 1 \langle
\alpha \rangle \amalg \beta \rangle$. The complex scheme of $A$ is
encoded by $\langle \mathcal{J} \amalg 1_{\epsilon} \langle \alpha_+ \amalg
\alpha_- \rangle \amalg \beta_+ \amalg \beta_- \rangle$ where
$\epsilon \in \{ +,- \}$ is the sign of the non-empty oval;
$\alpha_+, \alpha_-$ are the numbers of positive and negative ovals
among the $\alpha$; $\beta_+, \beta_-$ are the numbers of positive
and negative ovals among the $\beta$ (remember that all signs are defined
with respect to the orientation of $\mathcal{J}$).
\begin{description}
\item[Rokhlin-Mishachev formula:]
{\em If\/} $m = 2k + 1$, {\em then\/}
\begin{displaymath}
2(\Pi_+ - \Pi_-) + (\Lambda_+ - \Lambda_-) = L - 1 - k(k+1)
\end{displaymath}
\item[Fiedler theorem:]
{\em Let $\mathcal{L}_t = \{L_t, t \in [0,1]\}$ be a pencil of real lines
based in a point $P$ of $\mathbb{R}P^2$. Consider two lines $L_{t_1}$
and $L_{t_2}$ of $\mathcal{L}_t$, which are tangent to $\mathbb{R}A$ at
two points $P_1$ and $P_2$, such that $P_1$ and $P_2$ are related by a pair
of conjugated imaginary arcs in $\mathbb{C}A \cap (\bigcup L_t)$.

Orient $L_{t_1}$ coherently to $\mathbb{R}A$ in $P_1$, and transport
this orientation through $\mathcal{L}_t$ to $L_{t_2}$.
Then this orientation of $L_{t_2}$ is compatible to that of $\mathbb{R}A$
in $P_2$.\/}
\end{description}

A sequence of ovals that are connected one to the next by pairs of conjugated imaginary arcs will be called a
{\em Fiedler chain\/}. Such chains are easily spotted when the pencil of lines is {\em maximal\/}, that is to say has
alternatively $m$ and $m-2$ real intersection points with $\mathbb{R}A$. See \cite{or1} for the close connection of
Fiedler's theorem with Orevkov's quasi-positive braids. 

\subsection{Results}
After systematic constructions, Korchagin  \cite{ko3} proposed several conjectures about the $M$-curves of degree $9$, 
here is the one dealing with the $M$-curves with three nests:

\begin{conjecture}
{\em Let $C_9$ be an $M$-curve of degree 9 with real scheme
$\langle \mathcal{J} \amalg 1 \langle \alpha_1 \rangle
\amalg 1 \langle \alpha_2 \rangle \amalg 1 \langle \alpha_3 \rangle \amalg \beta\rangle$. 
At least two of the $\alpha_i, i = 1, 2, 3$ are odd.\/}
\end{conjecture}

Let $O_1, O_2, O_3$ be the non-empty ovals.
We call {\em nest\/} $\mathcal{O}_i$ each configuration  $1 \langle \alpha_i \rangle$
formed by $O_i$ and its interior ovals. Unfortunately, the word nest was used
with different meanings by the several authors in the past. Take good care that the definition used here differs from the more 
standard one from \cite{or1}, \cite{or6} or \cite{flt3}.

In the present paper, we prove a part of the conjecture.

\begin{description}
\item[Theorem 1] 
{\em Let $C_9$ be an $M$-curve of degree 9 with real scheme
$\langle \mathcal{J} \amalg 1 \langle \alpha_1 \rangle
\amalg 1 \langle \alpha_2 \rangle \amalg 1 \langle \alpha_3 \rangle \amalg \beta\rangle$. 
At least one of the $\alpha_i, i = 1, 2, 3$ is odd.\/}
\end{description}

Theorem 1 excludes the 53 real schemes $\langle \mathcal{J} \amalg 1 \langle \alpha_1 \rangle
\amalg 1 \langle \alpha_2 \rangle \amalg 1 \langle \alpha_3 \rangle \amalg \beta\rangle$
($\alpha_1 + \alpha_2 + \alpha_3 + \beta = 25$, $\alpha_1 \leq \alpha_2 \leq \alpha_3$)
with $\alpha_1, \alpha_2, \alpha_3$ even. Among them, the 12 ones with 
$\beta = 1$ had already been excluded by A. Korchagin in \cite{ko6}.
Our proof involves some classical tools: 
Bezout's theorem with auxiliary lines and conics, the Rokhlin-Mishachev formulas, Fiedler's theorem.
We use supplementarily rational cubics, quartics (single curves or pencils of such curves), and Orevkov's complex orientation
formulas from \cite{or6}. All of the arguments, and hence the statements, 
are also valuable for pseudo-holomorphic curves.

We prove also a few results on complex orientations and rigid isotopy for the $M$-curves $C_9$ with three nests.

\section{First properties}
\subsection{Descriptive lemmas and definitions}
Let $C_9$ be an $M$-curve of degree 9.
Given an empty oval $X$ of $C_9$, we often will have to consider one point
chosen in the interior of $X$. For simplicity, we shall call this point
also $X$. In the following, it will be clear from the context whether we speak of the oval
or of the point $X$.
We denote the pencil of lines based in $X$ by $\mathcal{F}_X$.
Let $[XY]$, and $[XY]'$ be the two segments of line determined by $X$ and $Y$, cutting
$\mathcal{J}$ respectively an even and an odd number of times. We say that $[XY]$ is the
{\em principal segment\/} determined by $X, Y$. Let $X, Y, Z$ be three ovals of $C_9$. Corresponding three
points $X, Y$ and $Z$ determine four triangles of $\mathbb{R}P^2$.
We will call {\em principal triangle\/} and denote by $XYZ$ the triangle
whose edges are the principal segments $[XY]$, $[YZ]$ and $[ZX]$.
Let $C_2$ be a conic passing through five points $A, B, C, D, E$ in this ordering. We write $C_2 = ABCDE$. 
If $F$ lies in the interior of $C_2$, we write $F < C_2$, otherwise we write $F > C_2$.

\begin{definition}
An ordered group of empty ovals $F_1, \ldots, F_n$ of $C_9$
lies in a {\em convex position\/} if for each triple
$F_i, F_j, F_k$, the principal triangle $F_iF_jF_k$ does not contain
any other oval of the group and $F_1, \ldots, F_n$ are the successive
vertices of $\bigcup F_iF_jF_k$ (the {\em convex hull of the group\/}).
\end{definition}
\begin{definition}
Let $O$ be a non-empty oval of $C_9$.
We say that $C_9$ has a {\em jump in $O$ determined by $A, B, C, D$\/} if there exist two empty
ovals $B$ and $C$ inside of $O$, and two empty ovals $A$ and $D$
outside of $O$, such that: $A$ lies inside of an oval $O'$ different
from $O$, and a line passing through $A$ and $D$ separates $B$
and $C$ in Int($O$), see Figure 1.
\end{definition}

\begin{figure}
\centering \psfig{file=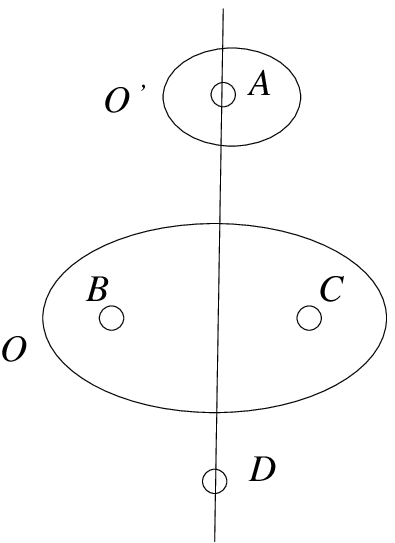}
\caption{jump}
\end{figure}

Let $C_9$ have real scheme
$\langle \mathcal{J} \amalg 1 \langle \alpha_1 \rangle \amalg 1 \langle \alpha_2 \rangle \amalg 
1 \langle \alpha_3 \rangle \amalg \beta \rangle$
and let $A_i, i=1, 2, 3$ be empty ovals of $\mathcal{O}_i$. The lines $(A_1A_2)$, $(A_2A_3)$, $(A_3A_1)$, and the pseudo-line $\mathcal{J}$
determine four triangles $T_0, T_1, T_2, T_3$ and three quadrangles $Q_1, Q_2, Q_3$ in
$\mathbb{R}P^2$.
Note that, by Bezout's theorem, $\mathcal{J}$ does not cut $T_0$, see Figure 2.

\begin{figure}
\centering \psfig{file=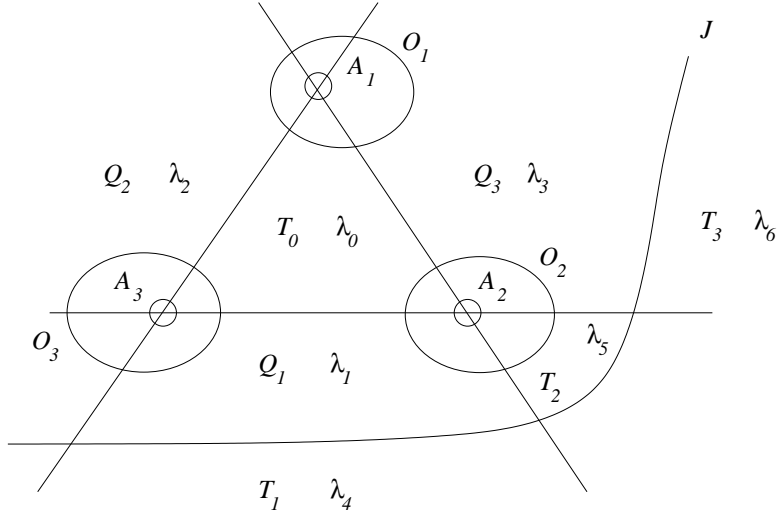}
\caption{Curve $C_9$ with three nests}
\end{figure}

The Lemmas 1, 3, 4, 7 hereafter are proven in \cite{flt2}.  

\begin{lemma}
If $C_9$ has a jump in $O_i$ determined by four ovals $A, B, C, D$, then $D$ is exterior.
\end{lemma}

\begin{lemma} 
Let $\{i, j, k\} = \{1, 2, 3\}$. All of the lines through two ovals interior to $O_i$ cut the same segment of 
line $[A_jA_k]$ or $[A_jA_k]'$.
\end{lemma}

{\em Proof\/} Let $A, B, C$ be three ovals in $Int(O_i)$. By Bezout's theorem
with the conic through $A, B, C, A_j, A_k$, the lines $(AB), (AC), (BC)$ must all cut the same segment
$[A_jA_k]$ or $[A_jA_k]'$. $\Box$

\begin{definition}
A non-empty oval $O_i$ of $C_9$ is {\em separating\/} if any line $(AA')$ 
joining two ovals of $Int(O_i)$ cuts the principal segment $[A_jA_k]$.
Otherwise, $O_i$ is {\em non-separating\/}
\end{definition}

\begin{lemma}
Let $C_9$ have a jump, say in $O_3$, determined by $A_1, B, C, D$. Then, 
up to permutation of $B, C$, the ovals $A_1, A_2, C, D, B$ lie in convex 
position. 
\end{lemma}

\begin{lemma}
Assume $C_9$ has a jump in $O_3$ determined by $A, B, C, D$, with $A$  
interior to $O_i, i=1$ or $2$. Let $A'$ be any other empty oval in $\mathcal{O}_1 \cup \mathcal{O}_2$.
Then, $A', B, C, D$ also give rise to a jump.
\end{lemma}

As an immediate consequence of Lemma 3, we get:

\begin{lemma}
If $C_9$ has a jump in $O_3$, then $O_3$ is non-separating.
\end{lemma}

\begin{lemma}
The pencils of lines based in empty ovals of $\mathcal{O}_i \cup \mathcal{O}_j$ sweeping out $O_k$
give all rise to the same Fiedler chain of empty ovals, the {\em chain of $O_k$\/}.
\end{lemma}

{\em Proof\/} 
Assume $k = 3$ and let $A, A'$ be empty ovals in $\mathcal{O}_1 \cup \mathcal{O}_2$.
Consider the pencils of lines $\mathcal{F}_{A}$ and $\mathcal{F}_{A'}$ sweeping out $O_3$, and   
three empty ovals $D, E, F$ in $\mathcal{O}_3$ met successively by the pencil based in $A$. 
By Lemmas 1-2, $A'$ must lie in the same triangle $DEF$ as $A$. 
Thus the interior ovals are met with the same ordering by the pencils $\mathcal{F}_A$ and $\mathcal{F}_{A'}$.
Let $O_3$ have a jump. By Lemma 4, the pencils $\mathcal{F}_A$ and $\mathcal{F}_{A'}$ sweep out the same sets of ovals. 
If $D_1, D_2$ are two of these ovals, then the line $(D_1D_2)$ cuts the non-principal segment $[AA']'$.
This follows from Lemmas 3-5 if at least one of the $D_i$ is interior.
Let $D_1, D_2$ be exterior. Assume that $\mathcal{F}_A$ sweeps out successively $B, D_1, D_2, C$
and $\mathcal{F}_{A'}$ sweeps out successively $B, D_2, D_1, C$. Then, the conic $ABA'D_2D_1$ cuts $C_9$ at
20 points, which is a contradiction.   
Therefore, both pencils sweep out successively the same sequence of ovals $B, D_1, D_2, C$ or $B, D_2, D_1, C$.
$\Box$

\begin{definition}
Let $O$ be a non-empty oval of $C_9$ and $S$ be an oval of $C_9$ lying 
inside of another non-empty oval $O'$. 
The curve $C_9$ has {\em $n$ jumps in $O$ with repartition\/}
($l_1, \ldots, l_{2n+1}$) if a pencil of lines $\mathcal{F}_S$ sweeping
out $O$ meets successively $2n + 1$ groups of ovals, which are
situated alternatively in, out, \dots, in Int($O$) and have cardinals
$l_1, \ldots, l_{2n+1}$.
\end{definition}

It follows from Lemma~6 that the number of jumps in $O$ and their
repartition does not depend on the choice of $S$. Thus Definition~4 is
correct.

\begin{lemma}
The curve $C_9$ has at most one jump. 
\end{lemma}

If $C_9$ has a jump, we fix the convention that it is in $O_3$.   

\begin{lemma}
Let $C_9$ be as in Figure 2, and choose a supplementary empty oval $A'_i$ in one of the nests
$\mathcal{O}_i$.
The lines $(A'_iA_j)$, $(A'_iA_k)$, $(A_jA_k)$, and $\mathcal{J}$ give rise to new triangles $T'_0$, $T'_1$, $T'_2$, $T'_3$.
Let $E$ be another empty oval.
\begin{enumerate}
\item
If $O_i$ is non-separating, then $E \in T_l \iff E \in T'_l$ for $l = 0, \dots 3$.
\item
If $O_i$ is separating, then: $E \in T_0 \cup T_i \iff E \in T'_0 \cup T'_i$; and  
$E \in T_l \iff E \in T'_l$ for $l = j, k$. If $E \in T_i \cap T'_0$, then $E$ lies inside of $O_i$.
\end{enumerate}
\end{lemma}

The proof is straightforward, applying Bezout's theorem with the conic through $A_i, A'_i, A_j, A_k, E$. $\Box$
  
Let $C_9$ have a jump in $O_3$ arising from $B, C, D$, and let $A_1 \in \mathcal{O}_1,$ 
$A_2 \in \mathcal{O}_2$ be the first and the last empty oval of $\mathcal{O}_1 \cup \mathcal{O}_2$ met
by the pencil of lines $\mathcal{F}_D: B \to O_1 \to O_2 \to C$. The sector $(DA_1, DA_2)$ swept out by this pencil
is the union of two triangles $DA_1A_2$, one of them being principal. By Bezout's Theorem, $O_3$ can intersect only
one of these triangles. 

\begin{definition}
Let $C_9$ have a jump in $O_3$ arising from $B, C, D$.
If $O_3$ intersects the principal triangle $DA_1A_2$, then we say that $D$ is {\em front\/}; 
otherwise we say that $D$ is {\em back\/}.
If all of the exterior ovals in the chain of $O_3$ are front, we say that $O_3$ is {\em crossing\/}; 
if all of the exterior ovals in the chain of $O_3$ are back, we say that $O_3$ is {\em non-crossing\/}.
\end{definition}

\begin{lemma}
Let $C_9$ have a jump in $O_3$ arising from empty ovals $A_3 = B, C, D$, such that $A_1, A_2, C, D, B$ lie in convex
position (see Figure 3). 
\begin{enumerate}
\item
If $D$ is front, there are no ovals in $T_3$.
\item
If $D$ is back, there are no ovals in $T_0 \cup T_1 \cup T_2$. 
\end{enumerate}
\end{lemma}

{\em Proof\/} 
Let $C_2$ be the conic through $A_1, A_2, B, C$ and another empty oval $E$. 
The triple $(A_1, A_2, B)$ determines triangles $T_l, l = 0, \dots 3$,
and the triple $(A_1, A_2, C)$ determines triangles $T'_l, l = 0, \dots 3$.  
Recall that by Lemma 5, $O_3$ is non-separating.
\begin{enumerate}
\item
Let $E \in T_3$. By Lemma 8(1) with $i = 3$, one has $E \in T_3 \cap T'_3$, thus $C_2 = A_1EA_2CB$.
By Bezout's theorem with $C_9$, the arc $CB$ of $C_2$ lies inside of $O_3$.
Thus, $D > C_2$, and $E < A_1A_2CDB$. The conic $A_1A_2CDB$ cuts $C_9$ at 20 points. This is a contradiction.  
\item
Let $E \in T_0$, then, again by Lemma 8(1) with $i = 3$, one has $E \in T_0 \cap T'_0$. 
The conic $C_2 = A_1EA_2BC$ cuts $C_9$ at 20 points.  
Let $E \in T_1 \cup T_2$. By symmetry, we may suppose that $E \in T_1$.
As $E \in T_1 \cap T'_1$, one has a priori $C_2 =A_1A_2ECB$ or $A_1A_2BCE$. 
By Bezout's theorem with $C_9$, one must have $C_2 = A_1A_2ECB$, and the arc $CB$ of $C_2$ lies inside of $O_3$.
Thus $D < C_2$ and $E < A_1A_2CDB$. The conic $A_1A_2CDB$ cuts $C_9$ at 20 points. This is a contradiction.
$\Box$
\end{enumerate}

\begin{figure}
\centering \psfig{file=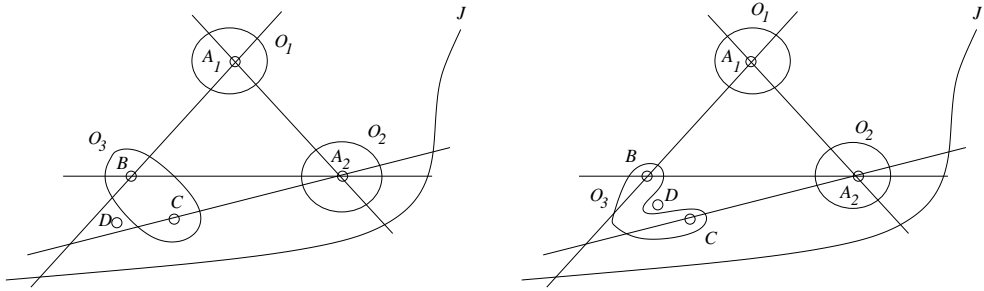}
\caption{left: $D$ is front, right: $D$ is back}
\end{figure}

\subsection{Complex orientions}

Let $a_i^\pm$ be the numbers of positive and negative interior ovals of the
nest $\mathcal{O}_i$. Let $A$ be any empty oval of $\mathcal{O}_j \cup \mathcal{O}_k$.
Recall that the chain of $O_i$ is the Fiedler chain of empty ovals
arising from the pencil $\mathcal{F}_{A}$ sweeping out $\mathcal{O}_i$,
this chain is independant of the choice of $A$ by Lemma 6.
There is at most one jump in $O_i$, thus $\vert a_i^+ - a_i^- \vert \leq 2$. The equality is achieved if and only if
$O_i$ has a jump with repartition $l_1, l_2, l_3$, with each $l_n, n=1, 2, 3$
odd.
Let us call {\em principal ovals\/} the ovals $O_1, O_2, O_3, A_1, A_2, A_3$, and
{\em base ovals\/} the empty principal ovals $A_1, A_2, A_3$.
Let $\epsilon_n, n=1, 2, 3, 4, 5, 6$, $\epsilon_n \in \pm 1$ be 
the respective contributions of these 6 ovals to $\Lambda_+ - \Lambda_-$.
Let $\lambda_0, \lambda_1, \lambda_2, \lambda_3, \lambda_4, \lambda_5, 
\lambda_6$ be the contributions to $\Lambda_+ - \Lambda_-$ brought
respectively by the non-principal ovals of the zones $T_0$, $Q_1, Q_2, Q_3,
T_1, T_2, T_3$. 

\begin{lemma}
One has:
\begin{displaymath}
\lambda_0 + \lambda_1 - \lambda_4 = -\frac{1}{2}
(\epsilon_3 + \epsilon_6 + 
\epsilon_2 + \epsilon_5)
\end{displaymath}
\begin{displaymath}
\lambda_0 + \lambda_2 - \lambda_5 = -\frac{1}{2}
(\epsilon_3 + \epsilon_6 + 
\epsilon_1 + \epsilon_4)
\end{displaymath}
\begin{displaymath}
\lambda_0 + \lambda_3 - \lambda_6 = -\frac{1}{2}
(\epsilon_2 + \epsilon_5 + 
\epsilon_1 + \epsilon_4)
\end{displaymath}
\begin{displaymath}
3\lambda_0 + \lambda_1 + \lambda_2 + \lambda_3 -  \lambda_4 - \lambda_5 - 
\lambda_6 + \sum{\epsilon_i} = 0
\end{displaymath}
\begin{displaymath}
\lambda_0 -\lambda_4 -\lambda_5 -\lambda_6 = 
-\frac{1}{2}( \Lambda_+ - \Lambda_-) = \Pi_+ - \Pi_- - 4
\end{displaymath}
\end{lemma}

{\em Proof\/} Apply Fiedler's Theorem to the pencils of lines $\mathcal{F}_{A_1}: A_3 \to T_0 \cup Q_1 \cup T_1 \to A_2$,
$\mathcal{F}_{A_2}: A_1 \to T_0 \cup Q_2 \cup T_2 \to A_3$ and  $\mathcal{F}_{A_3}: A_2 \to T_0 \cup Q_3 \cup T_3 \to A_1$. 
Subtracting $\lambda_0 + \lambda_1 + \lambda_2 + \lambda_3 +  \lambda_4 + \lambda_5 + 
\lambda_6 + \sum{\epsilon_i} = \Lambda_+ - \Lambda_-$ from the fourth identity, and combining with Rokhlin-Mishachev's formula
yields the last identity.
$\Box$

\section{Cremona transformation}

Let $C_9$ be, as in the previous section, an $M$-curve of degree 9 with real scheme
$\langle \mathcal{J} \amalg 1 \langle \alpha_1 \rangle
\amalg 1 \langle \alpha_2 \rangle \amalg 1 \langle \alpha_3 \rangle \amalg \beta \rangle$.
Let us perform a Cremona transformation $cr: (x_0: x_1: x_2) \to (x_1x_2: x_0x_2: x_0x_1)$ with base points 
$A_1, A_2, A_3$. We shall denote the respective images of the lines $(A_1A_2), (A_2A_3),(A_3A_1)$ by $A_3, A_1, A_2$.
For the other points, we use the same notation as before $cr$.
The curve $C_9$ is mapped onto a curve $C_{18}$ of degree 18 with three singular points.
We shall call {\em main part of\/} $C_{18}$ the piece formed by the images of 
$\mathcal{J}$ and the principal ovals, see Figure 4 where $cr(A_i)$ and $cr(O_i)$ stand for the images of the
{\em ovals\/} $A_i$ and $O_i$. 

An oval $A$ of $C_{18}$ will be said to be {\em interior, exterior, positive\/} or {\em negative\/} if its
preimage is. Let $\mathcal{O} = cr(\mathcal{J})$. One has $Int(\mathcal{O}) = cr(  T_0 \cup T_1 \cup T_2 \cup T_3)$, 
$Ext(\mathcal{O}) = cr( Q_1 \cup Q_2 \cup Q_3)$. The ovals of $Int(\mathcal{O})$ and their preimages 
will be called {\em triangular ovals\/}; the ovals of $Ext(\mathcal{O})$
and their preimages will be called {\em quadrangular ovals\/}. 
A configuration of ovals of $C_{18}$ lies in {\em convex position\/} if there exists a line $L$ such that the
ovals lie in convex position in the affine plane $\mathbb{R}P^2 \setminus L$.
Consider the injective pairs formed by the triangular ovals and $\mathcal{O}$. One has
$\pi_+ - \pi_-(\mathcal{O}) = \lambda_4 + \lambda_5 + \lambda_6 - \lambda_0$.

\begin{figure}
\centering \psfig{file=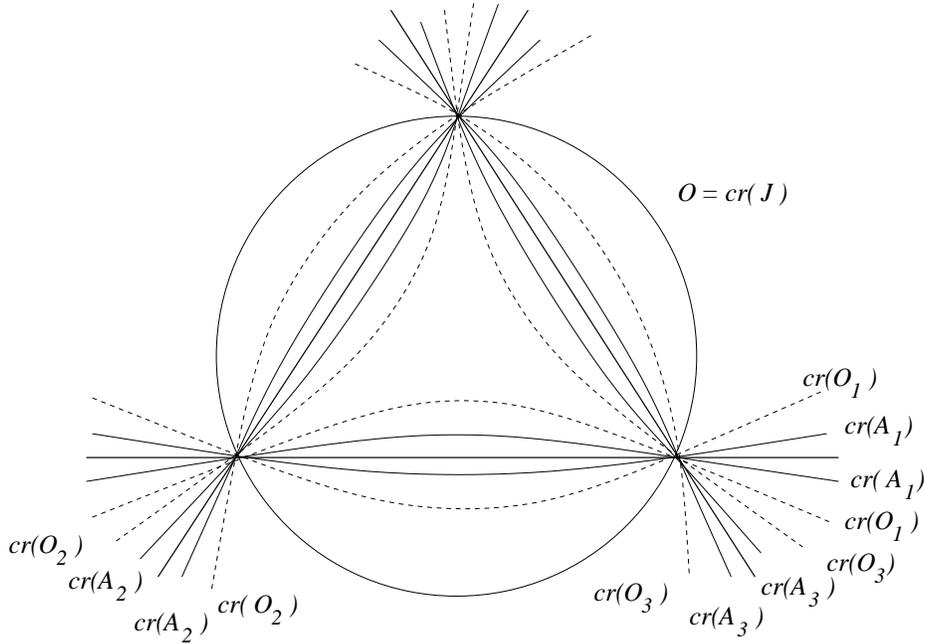}
\caption{Main part of $C_{18}$}
\end{figure}

\subsection{Lemmas using auxiliary conics}

\begin{definition}
Let $\{i, j, k\} = \{1, 2, 3\}$.
A base line $(A_iA_j)$ and a conic $C_2$ are mutually {\em maximal\/} if $C_2$ cuts $(A_iA_j)$ and $C_2$ cuts 
each component $cr(O_k)$ and $cr(A_k)$ at four points.
\end{definition}

\begin{lemma}
If a base point $A_i$ lies inside of a conic $C_2$, then the two base lines $(A_iA_j)$ and $(A_iA_k)$ are maximal with respect 
to $C_2$.
\end{lemma}

{\em Proof\/}
Let $P_1, Q_1 = C_2 \cap (A_2A_3)$,  $P_2, Q_2 = C_2 \cap (A_1A_3)$ and  $P_3, Q_3 = C_2 \cap (A_1A_2)$
If $A_i$ lies inside of $C_2$, then $C_2$ meets $P_k, P_j, Q_k, Q_j$ in this ordering. Each arc joining two consecutive
points cuts $cr(O_k), cr(A_k), cr(O_j)$ and $cr(A_j)$. $\Box$

\begin{lemma}
Let $C_2$ be a conic passing through five ovals $E_1, \dots E_5$ of $C_{18}$. Let $A_i, A_j$ be two of the base points, lying outside of 
$C_2$, such that the line $(A_iA_j)$ cuts $C_2$.
If any of the following conditions is verified, then $(A_iA_j)$ is maximal with respect to $C_2$.
\begin{enumerate}
\item
Each arc of $C_2 \setminus (C_2 \cap A_iA_j)$ passes through an oval $E_m$ that is exterior to $cr(O_k)$, or
cuts $\mathcal{O}$.
\item
$Int(C_2) \cup Int(\mathcal{O})$ is orientable and $A_i, A_j$ lie on different arcs of $\mathcal{O} \setminus (\mathcal{O} \cap C_2)$,
\end{enumerate}
\end{lemma}

We leave the proof to the reader.


\begin{lemma}
Let $C_2$ be a conic passing through five ovals $E_1, \dots E_5$ of $C_{18}$, and having at least four intersection points
with $\mathcal{O}$. Then one of the three base lines, say $A_1A_3$ is non-maximal with respect to $C_2$,
and the points $A_1, A_3$ lie outside of $C_2$.
\end{lemma}

{\em Proof\/} 
If the three base lines are maximal with respect to $C_2$, then $C_2$ cuts the images of the principal ovals at
$24$ points, $\mathcal{O}$ at four points, and the union $\cup E_i, i=1, \dots, 5$ at $10$ points.
This is a contradiction. A base line, say $A_1A_3$ is non-maximal; the points $A_1, A_3$ must lie outside of $C_2$ by Lemma 11. $\Box$

\begin{lemma}
Let $C_2 =  E_1E_2E_3E_4E_5$ be a conic (passing through none of the base points), satisfying: 
$Int(\mathcal{O}) \cup Int(C_2)$ is orientable, $E_1, E_3, E_5$ are quadrangular, 
$E_2, E_4$ are triangular. Consider the set $\mathcal{S}$ of arcs of
$\mathcal{O} \setminus (\mathcal{O} \cap C_2)$, they are isotopic with fixed endpoints to arcs of $C_2$.
\begin{enumerate}
\item
There exist four arcs $s_1, \dots, s_4$ of $\mathcal{S}$  
whose endpoints lie on the following arcs of $C_2$: $E_iE_{i+1}, E_{i+1}E_{i+2}$ for $i = 1, 2, 3$, 
and  $E_4E_5, E_1E_2$ for $i = 4$. The arcs $s_1, s_3$ are exterior and the arcs $s_2, s_4$ are interior to $C_2$. 
Consider any other arc $s'$ of $\mathcal{S}$. Both endpoints of $s'$ lie on the same arc
$E_1E_2$, $E_2E_3$, $E_3E_4$ or $E_4E_5$ of $C_2$.
\item
One of the base lines, say $A_1A_3$, is non-maximal for $C_2$ and $A_1, A_3$ lie both on the same arc $s = s_1$ or
$s_3$ of  $\mathcal{S}$
\end{enumerate}
\end{lemma}


{\em Proof\/}
Lemma 13 implies that a base line, say $(A_1A_3)$, is non-maximal for $C_2$; by Lemma 12 (2), the points $A_1, A_3$
lie on the same arc $s$ of $\mathcal{S}$, exterior to $C_2$.
By Bezout's theorem, the lines $(E_iE_j)$ have at most two intersection points with $\mathcal{O}$.
Each arc $E_iE_{i+1}$, $i=1, \dots 4$ of $C_2$ cuts $\mathcal{O}$ an odd number of times.
The existence of arcs $s_i, i = 1, 2, 3$ of $\mathcal{S}$ as required in 1 follows immediately from the hypotheses.
By Bezout's theorem with the lines $(E_iE_{i+1}), i = 1, \dots 4$, any other arc $\sigma$ of
$\mathcal{S}$ with one endpoint on the arc $E_2E_3$ or $E_3E_4$ of $C_2$ must have its second endpoint on the
same arc of $C_2$. Let $P$ be the endpoints of $s_1$ on the arc $E_1E_2$ of $C_2$, and let $Q$ be the endpoint of $s_3$ on the
arc $E_4E_5$ of $C_2$. One arc $PQ$ of $\mathcal{O}$ is a union $s_1 \cup s_2 \cup s_3 \cup \sigma_1 \cup \dots \cup \sigma_n$  of arcs 
of $\mathcal{S}$, where each of the $\sigma_k$ has both endpoints on the same arc $E_2E_3$ or $E_3E_4$ of $C_2$.
The other arc $PQ$ of $\mathcal{O}$ is a union $s_4 \cup \tau_1 \cup \dots \cup \tau_m$, such that: 
$\tau_1, \dots \tau_m \in \mathcal{S}$, and each $\tau_k$ has both endpoints on
the same arc $E_1E_2$ or $E_4E_5$ of $C_2$; $s_4$ has germs near the endpoints interior to $C_2$ and possibly supplementary pairs
of intersections with the arc $E_5E_1$ of $C_2$.
We will prove that the arc $s$ containing $A_1, A_3$ must be $s_1$ or $s_3$.
Assume that $s$ is another arc of $\mathcal{S}$. The endpoints of $s$ are on the same
arc of $C_2$, up to symmetry, we may assume that this arc is $E_1E_2$, $E_2E_3$ or $E_5E_1$.
Let us start with the case where the arc is $E_1E_2$. Consider the conics $C_2(A_1) = A_1E_2E_3E_4E_5$ and 
$C_2(A_3) =  A_3E_2E_3E_4E_5$, see Figure 5. One of the two base points lies in the interior of the conic passing through 
the other, say $A_3$ lies in the interior of $C_2(A_1)$; the other base point $A_2$ may be inside or outside. 
Then, $C_2(A_1)$ cuts $\mathcal{O}$ at 6 points, each of $cr(O_i)$, $cr(A_i), i = 1, 2, 3$ at 4 points, and each of the ovals
$E_2, \dots, E_5$ at 2 points. Hence in total 38 intersection points with $C_{18}$.
This is a contradiction.
The cases where $s$ has its endpoints on $E_2E_3$ or on $E_5E_1$ are similar, letting $C_2(A_i), i=1, 3$ be 
respectively: $A_iE_3E_4E_5E_2$ and $A_iE_1E_2E_3E_4$. To finish the proof, one needs to show that
the arc $E_5E_1$ of $C_2$ cannot cut $\mathcal{O}$. Assume the contrary. Up to symmetry, we may assume that $s = s_1$.
We find again a contradiction, using conics $C_2(A_i) = A_iE_3E_4E_5E_1$. 
$\Box$



\begin{figure}
\centering \psfig{file=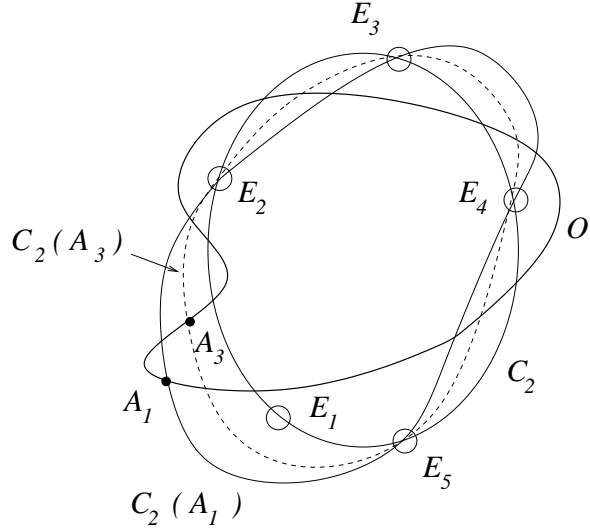}
\caption{$s$ with endpoints on the arc $E_1E_2$ of $C_2$} 
\end{figure}


\begin{lemma}
Let the curve $C_{18}$ contain a configuration of six ovals 
$B_1$, $D_3$, $B_2$, $D_1$, $B_3$, $D_2$ lying in convex position, such that
$B_1, B_2, B_3$ are triangular, and the segments $[B_iB_j]$ interior to the convex hull are also interior to
$\mathcal{O}$. Then, one of the ovals $D_1, D_2, D_3$ must also be triangular.
\end{lemma}


{\em Proof\/} 
Let $B_i, D_i, i = 1, 2, 3$ satisfy the conditions of the Lemma and assume that the $D_i$ are all
quadrangular. Denote by $\mathcal{H}$ the hexagon bounding the convex hull. 
Each of the six edges $[B_iD_j]$ of $\mathcal{H}$ cuts $\mathcal{O}$ once.
Let $P_1, Q_1, P_2, Q_2, P_3, Q_3$ be the successive intersections of $\mathcal{O}$ with the edges $[D_2B_1]$, $[B_1D_3], \dots
[B_3D_2]$. Connect pairwise these six points by arcs of $\mathcal{O}$ interior to $\mathcal{H}$, these arcs must be 
$Q_1P_2, Q_2P_3, Q_3P_1$.
Connect pairwise the six points with exterior arcs. By 
Bezout's theorem with the lines $(D_iD_j)$, these arcs must be $Q_1P_1, Q_2P_2, Q_3P_3$. 
Consider the three conics $B_1D_3D_1B_3D_2$, $B_3D_2D_3B_2D_1$ and $B_2D_1D_2B_1D_3$.  
Bezout's theorem with the six lines supporting $\mathcal{H}$ implies that $Int(\mathcal{O}) \cup Int(C_2)$ is orientable. 
By Lemma 14(1), the arc $D_kD_i$ of $B_iD_kD_iB_kD_j$ does not cut $\mathcal{O}$. Each oval $B_i$ lies inside of the conic 
$C_2$ determined by the other five ovals.  
Hence, each exterior arc of $\mathcal{O} \setminus (\mathcal{O} \cap C_2)$ is entirely contained in another conic.
Thus, any choice of base points leads to a contradiction with Lemma 14(2), see Figure 6. $\Box$ 

\subsection{Inequalities}

In the proofs of the next two statements, we consider conics passing through some empty ovals of $C_{18}$.
Several times, we find a conic that is maximal with respect to the three base lines. The maximality follows always 
from Lemma 12 (1): each base line separates on this conic a pair of exterior ovals.  
Let $L, L'$ be two lines and $D$ be a point, we denote by $(L, L', \hat D)$ the sector $(L, L')$ that does not contain $D$.

\begin{lemma}
Let the base ovals $A_1, A_2, A_3$ of $C_9$ be such that $T_0$ contains only exterior ovals
of $C_9$.
One has $\vert \lambda_0 \vert \leq 3$. If $\lambda_0 = \pm 3$, then the non-empty ovals are all
separating and $\sum \epsilon_i = \mp 6$. Moreover, one of the quadrangles is empty.
\end{lemma}

{\em Proof\/}
Note that the choice of $A_i$ is unique if $O_i$ is separating, and arbitrary if
$O_i$ is non-separating. 
Perform the Cremona transformation $cr$, and denote by $B_i, i=1, \dots, n$ the ovals of $C_{18}$ lying in $T_0$.
Assume there exist $B_i, B_j, B_k, B_l$ such that $B_l$ lies in the triangle $B_iB_jB_k$ that does not cut
the base lines. Consider the pencil of conics $\mathcal{F}_{B_iB_jB_kB_l}$,
the conics of this pencil are all maximal with respect to the three base lines.
They intersect $C_{18}$ at $36$ points; the other empty ovals of $C_{18}$ cannot be swept out. This is a contradiction.
Thus, the ovals of $T_0$ lie in convex position in this triangle, let $\Delta$ be their convex hull.
Let $D$ be any oval of $C_{18}$. Let $\mathcal{F}_D$ be: the complete pencil of lines based in $D$ if
$D$ is triangular; the portion of pencil based in $D$ and sweeping out $\mathcal{O}$ if $D$ is quadrangular.
Note that $\mathcal{F}_D$ is maximal, see the end of section 1.1.
Consider two successive vertices $B_i, B_j$ of $\Delta$ that are consecutive for some pencil $\mathcal{F}_{B_k}$.
Then, for any choice of another empty oval $D$, the ovals $B_i, B_j$ are consecutive for the pencil $\mathcal{F}_D$.
Indeed, assume there exists an oval $D$ that does not fullfill this condition. There exists a conic $C_2$ passing through 
$B_i, B_j, B_k, D$ and a fifth oval, that is maximal with respect to the three base lines and cuts $\mathcal{O}$
at four points, which is impossible.
Thus, the ovals of $T_0$, ordered by the convexity, split into a number $N$ of successive
Fiedler chains, whose base points need not to be specified.  
Assume $\vert \lambda_0 \vert \geq 3$. There exist three chains bringing the same contribution
$+1$ or $-1$ to $\lambda_0$. Choose $B_1, B_2, B_3$ to be three ovals distributed in the three chains. 
Let $[B_1B_2]$, $[B_1B_3]$ and $[B_2B_3]$ be the segments contained in $T_0$ determined by the three points.
For $\{ i, j, k \} = \{ 1, 2, 3 \}$, the pencil of lines $\mathcal{F}_{B_i}$ sweeping out $[B_jB_k]$ must meet an oval 
$D_i$ outside of $T_0$. One gets thus a configuration of six ovals. Consider the conics determined by the three ovals 
in $T_0$ and two of the other ovals.
By Bezout's theorem with $C_{18}$, these three conics are: $B_1D_3B_2D_1B_3$, $B_2D_1B_3D_2B_1$, $B_3D_2B_1D_3B_2$.
To each conic determined by five given points, we associate the convex pentagon inside the conic having these points as vertices.
Choose a line at infinity $L$ that does not cut any of the three pentagons. The points 
$B_1, D_3, B_2, D_1, B_3, D_2$ lie in convex position in the affine plane $\mathbb{R}P^2 \setminus L$. 
The hexagon $\mathcal{H} = B_1D_3B_2D_1B_3D_2$ gives rise to a natural cyclic ordering of the six lines supporting its edges.
Let $Z_l, l \in \{1, \dots, 6\}$ be the six triangles that are supported by triples of consecutive lines,
such that the intersection of $Z_l$ and $\mathcal{H}$ is a common edge.
The base lines are distributed in the three zones: $(B_2B_3, B_2D_1, \hat B_1) \cup (B_3B_2, B_3D_1, \hat B_1)$,
$(B_3B_1, B_3D_2, \hat B_2) \cup (B_1B_3, B_1D_2, \hat B_2)$ and $(B_1B_2, B_1D_3, \hat B_3) \cup (B_2B_1, B_2D_3, \hat B_3)$.
(See the notation $(L, L', \hat D)$ introduced just before the statement of the Lemma.) 
Each triangle $Z_l$ has a nonempty intersection with $T_0$ and
cuts only one base line: $Z_1, Z_2$ cut $(A_1A_2)$; $Z_3, Z_4$ cut $(A_2A_3)$; $Z_5, Z_6$ cut $(A_3A_1)$.
Thus, $Z_1 \cup Z_2 \subset T_0 \cup cr(Q_3 \cup T_3)$, $Z_3 \cup Z_4 \subset T_0 \cup cr(Q_1 \cup T_1)$, 
$Z_5 \cup Z_6 \subset T_0 \cup cr(Q_2 \cup T_2)$. 
Consider the pencil of conics $\mathcal{F}_{D_kB_jB_kB_i}$. By Bezout's theorem with $C_{18}$, the remaining empty ovals must be 
swept out in the portion $D_kB_i \cup B_jB_k \to D_kB_j \cup B_iB_k$. Thus, all of these ovals lie in
$\cup Z_l$. Let $E, F$ be two ovals in the same triangle with vertices $B_i$ and $D_j$. If the line $EF$ cuts 
one of the affine segments $[B_iB_j]$ or $[B_jB_k]$, then the conic $C_2$ determined by $E, F, B_1, B_2, B_3$ intersects $C_{18}$ at
$38$ points, which is a contradiction. 
Thus, there is a natural cyclic ordering of the empty ovals of
$C_{18}$ given by the pencils of lines $\mathcal{F}_{B_i}, i = 1, 2, 3$ sweeping out the six triangles (and the pencils
$\mathcal{F}_{D_i}, i = 1, 2, 3$ sweeping out the four triangles that do not have $D_i$ as vertex).
For $\{ i, j, k \} = \{ 1, 2, 3 \}$, the zone $cr(Int(O_i) \cap (Q_j \cup Q_k))$ does not cut $\cup Z_l$.
Thus, the interior ovals are all triangular. Otherwise stated, the non-empty ovals are all
separating. 
The $22$ ovals in the affine zone $\cup Z_l$ form a cyclic Fiedler chain splitting into 6 consecutive subchains in the zones
$T_0$, $cr(T_1 \cup Q_1)$,  $T_0$, $cr(Q_2 \cup T_2)$, $T_0$, $cr(Q_3 \cup T_3)$.
In this cyclic chain, an oval is say {\em positive\/} if and only if it is the image by $cr$ of a positive oval in $T_1 \cup T_2 \cup T_3$
or of a negative oval in $T_0 \cup Q_1 \cup Q_2 \cup Q_3$. 
Thus, $\lambda_0 + \lambda_1 + \lambda_2 + \lambda_3 - \lambda_4 - \lambda_5 - \lambda_6 = 0$.
Combining this with the fourth identity in Lemma 10, one gets: $2\lambda_0 = -\sum \epsilon_i$.
Any 6 consecutive ovals distributed in the 6 subchains lie in convex position. By Lemma 15, one of the zones $cr(T_i \cup Q_i)$
contains only triangular ovals. $\Box$


\begin{figure}
\centering \psfig{file=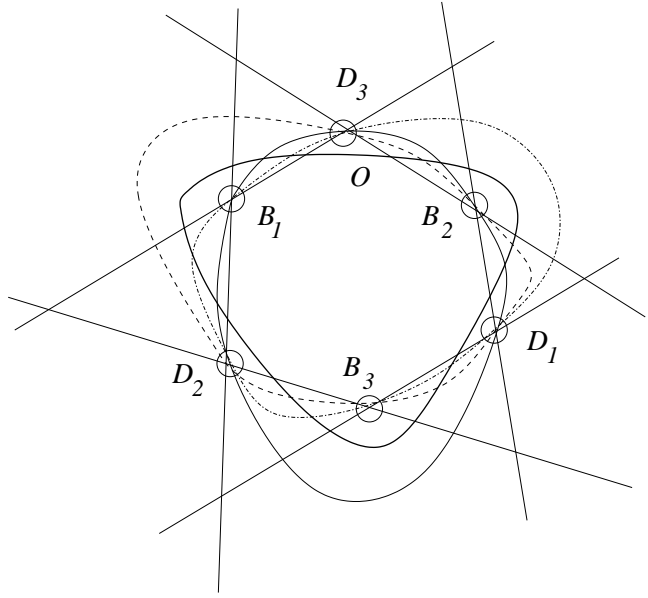}
\caption{The three conics}
\end{figure}

\begin{proposition}
Let $i \in \{1, 2, 3\}$, such that $T_i$ contains only exterior ovals. 
One has $\vert \lambda_{i+3} \vert \leq 3$. If $\vert \lambda_{i+3} \vert = 3$, then $\lambda_{i+3} = +3$ and
$\lambda_0 - \lambda_4 - \lambda_5 - \lambda_6 = -2$.
\end{proposition}

{\em Proof\/} 
Recall that for $i \in \{1, 2, 3\}$, the contribution of the ovals of $T_i$ to $\Lambda_+ - \Lambda_-$ is denoted by $\lambda_{i+3}$.
Perform the Cremona transformation $cr$, and denote by $B_1, \dots, B_n$ the exterior ovals of $C_{18}$ lying in $cr(T_i)$.
Let $B_q, B_r, B_s$ be three such ovals, note that the triangle $B_qB_rB_s$ whose edges do not cut the base lines is 
entirely contained in $cr(T_i)$. This triangle contains no other oval, by an argument similar to that used for Lemma 16. 
Thus, the ovals of $cr(T_i)$ lie in convex position in this zone, let $\Delta$ be their convex hull.
Consider two successive vertices $B_q, B_r$ of $\Delta$ that are consecutive for some pencil $\mathcal{F}_{B_s}$. 
Then, for any other choice of an empty oval $D$, the ovals $B_qB_r$ are consecutive for the maximal pencil $\mathcal{F}_D$. 
Thus, one may speak of a Fiedler chain of ovals in $cr(T_i)$, without refering to a base point.  
Assume that $\vert \lambda_{i+3} \vert \geq 3$, so there exist three Fiedler chains of ovals in $cr(T_i)$, bringing
each the same contribution $+1$ or $-1$ to  $\lambda_{i+3}$.
Choose $B_1, B_2, B_3$ to be three ovals distributed in the three chains. 
Let $[B_1B_2]$, $[B_1B_3]$ and $[B_2B_3]$ be the segments contained in $cr(T_i)$ determined by the three points.
For $\{q, r, s\} = \{1, 2, 3\}$, the pencil of lines $\mathcal{F}_{B_q}$ sweeping out $[B_rB_s]$ must meet an oval 
$D_q$ outside of $cr(T_i)$. One gets thus a configuration of six ovals. Consider the conics determined by 
the ovals $B_1, B_2, B_3$ and two of the other ovals. By Bezout's theorem with $C_{18}$, these three conics are: 
$B_1D_3B_2D_1B_3$, $B_2D_1B_3D_2B_1$, $B_3D_2B_1D_3B_2$.
Choose a line at infinity $L$ that does not cut any of the three interior pentagons. The points $B_1, D_3, B_2, D_1, B_3, D_2$
lie in convex position in the affine plane $\mathbb{R}P^2 \setminus L$. 
The hexagon $\mathcal{H} = B_1D_3B_2D_1B_3D_2$ gives rise to a natural cyclic ordering of the six lines supporting its edges.
Let $Z_l, l \in \{1, \dots, 6\}$ be the six triangles that are supported by triples of consecutive lines,
such that $Z_1 \cap \mathcal{H} = [B_1D_3]$, $Z_2 \cap \mathcal{H} = [D_3B_2]  \dots  Z_6 \cap \mathcal{H} = [D_2B_1]$. 
Similar arguments as in the proof of Lemma 16 yield the following statements:
all of the remaining ovals of $C_{18}$ lie in $\cup Z_l$. There is a natural cyclic ordering of the empty ovals of
$C_{18}$ given by the pencils of lines $\mathcal{F}_{B_q}, q = 1, 2, 3$ sweeping out the six triangles (and the pencils 
$\mathcal{F}_{D_q}, q = 1, 2, 3$ sweeping out the four triangles that do not have $D_q$ as vertex).
By Lemma 15, one of the $D_q$, say $D_3$, is triangular. 
Let $\{ i, j, k \} = \{ 1, 2, 3 \}$.
The ovals $B_1, B_2$ are separated from $D_3$ in $Int(\mathcal{O})$ by $[A_jA_k]$. Therefore, the line $(A_jA_k)$ lies 
in the zone $(B_1B_2, B_1D_3, \hat B_3) \cup (B_2B_1, B_2D_3, \hat B_3)$. 
The zone $cr(Q_i \cup T_i)$ is a triangle $T$ with vertices $A_1, A_2, A_3$.
If $D_q, q \in \{1, 2\}$ lies in $T$, then $D_q$ lies in $cr(Q_i)$.
The boundary of $cr(T_i)$ is a  topological circle $[A_jA_k] \cup (\mathcal{O} \cap T)$.
The oval $B_3$ lies inside of this circle, and the $D_q$ lie outside. 
The edge $[D_qB_3], q \in \{ 1, 2 \}$ does not cut $[A_jA_k]$, so it cuts $\mathcal{O}$. 
Hence, $D_1$ and $D_2$ are quadrangular. 
We may assume that $B_1$ and $B_3$ are respectively the first and the last oval met by the pencil of lines
$\mathcal{F}_{D_2}$ sweeping out $\mathcal{O}$. 
By Lemma 15, $Z_1 \cup Z_2$ contains only triangular ovals.
The complete pencil $\mathcal{F}_{D_2}$ starting at $B_1$ and sweeping out successively $Z_1, Z_2, Z_3, Z_4, Z_5 \cup Z_6$
meets four consecutive chains of ovals, alternatively triangular and quadrangular. 
Indeed, assume there exist more triangular and quadrangular subchains. 
There exist two ovals $E, F$ met successively by $\mathcal{F}_{D_2}$, such that:
$E, F \in Z_3$, $E$ is quadrangular, $F$ is triangular; or $E, F \in Z_4$,
$F$ is quadrangular, $E$ is triangular.
Let $E, F \in Z_3$. 
If the line $(EF)$ cuts the edge $[D_1B_3]$ of $\mathcal{H}$, then Bezout's theorem with the conic $EFD_1B_2B_3$ yields a contradiction.
Otherwise, $B_2, E, F, D_1, B_3, D_2$ lie in convex position and contradict Lemma 15.
The case where $E, F \in Z_4$ is similar.
The chain of triangular ovals containing $B_3$ lies in $Z_4 \cup Z_5$. As the line $(A_jA_k)$ does not cut these triangles, 
all of the ovals of this chain must lie in $cr(T_i)$.
The other chain of triangular ovals splits into three consecutive subchains containing respectively $B_1, D_3$ and $B_2$, such that
the first and the third chain lie in $cr(T_i)$, the second chains lies in $cr(T_0 \cup T_j \cup T_k)$.
Indeed, assume there exist more subchains in each zone. There exist two ovals $E, F$ met successively by $\mathcal{F}_{D_2}$
such that:
$E, F \in Z_1$, $E \in cr(T_0 \cup T_j \cup T_k)$, $F \in cr(T_i)$; or 
$E, F \in Z_2$, $F \in cr(T_0 \cup T_j \cup T_k)$, $E \in cr(T_i)$.
Let $E, F \in Z_1$. The line $(EF)$ cuts the edge $[D_3B_2]$ of $\mathcal{H}$ and the conic $EFD_3B_1B_2$ intersects $C_{18}$
at $38$ points, which is a contradiction. The case where $E, F \in Z_4$ is similar.
Thus  $\vert \lambda_{i+3} \vert \leq3$. 
If $\lambda_{i+3} = \pm 3$ then $\lambda_0 - \lambda_4 - \lambda_5 - \lambda_6 = \mp 2$.
On the other hand, the last identity of Lemma 10 implies that 
$\lambda_0 - \lambda_4 - \lambda_5 - \lambda_6 \leq 0$. This finishes the proof.
$\Box$ 


\section{$M$-curves with three nests and a jump}

Let $C_9$ be an $M$-curve of degree 9 with three nests and a jump arising from $B, C, D$, see Figure 3 where one 
makes $A_3 = B$. Let us assume that $C_9$ has some triangular ovals.
The 21 ovals of $C_9 \setminus \{O_1, O_2, O_3, A_1, A_2, A_3, C\}$
are swept out by the pencil of conics $\mathcal{F}_{A_1A_2A_3C}: A_1A_3 \cup A_2C \to A_1A_2 \cup A_3C \to A_2A_3 \cup A_1C$
if $O_3$ is crossing, and by the pencil of conics $\mathcal{F}_{A_1A_2A_3C}: A_1A_3 \cup A_2C \to A_1A_2 \cup A_3C$
if $O_3$ is non-crossing. In both cases, the 21 ovals are distributed in two Fiedler chains. Denote by $P, Q$ 
the pair of starting points of the chains, where $P$ is a points of tangency with $O_3$ and $Q$ is a point of tangency 
with $\mathcal{J}$. Denote by $P', Q'$ the endpoints of the chains, where $P'$ is a points of tangency with $O_3$ 
and $Q'$ is a point of tangency with $\mathcal{J}$. The pair of Fiedler chains is then: ($P \to P'$, $Q \to Q'$) or 
($P \to Q'$, $Q \to P'$).
The pencil of conics is mapped by $cr$ onto the pencil of lines $\mathcal{F}_C$, see Figure 7.

\begin{lemma}
All of the ovals in $T_0 \cup T_1 \cup T_2 \cup T_3$ are consecutive for the maximal portion of  
$\mathcal{F}_{A_1A_2A_3C}$ formed by conics intersecting $\mathcal{J}$.
If $O_3$ is crossing, they form a Fiedler chain: {\em triangular ovals\/} $\to P'$.
Thus, $\lambda_0 -\lambda_4 -\lambda_5 = 0$ or $-\epsilon_3$. 
If $O_3$ is non-crossing, they form a Fiedler chain: $P \to$ {\em triangular ovals\/}. 
Thus, $\lambda_6 = 0$ or $-\epsilon_3$.
\end{lemma}

\begin{figure}
\centering \psfig{file=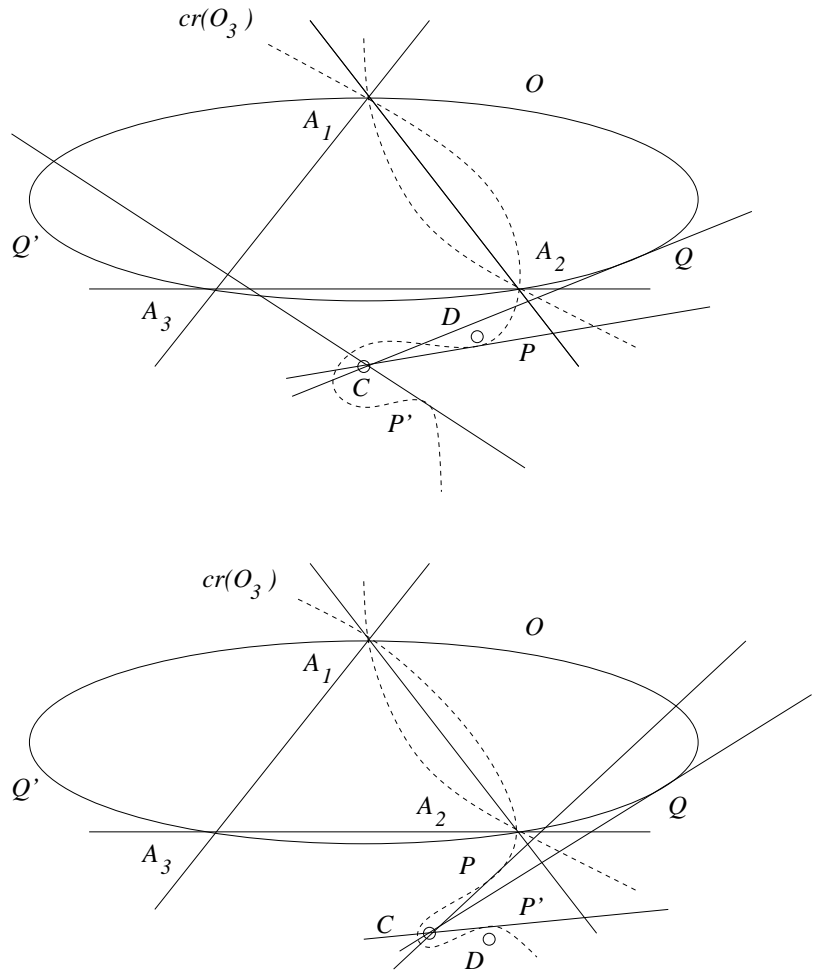}
\caption{Curve $C_{18}$ with the pencil of lines $\mathcal{F}_C$}
\end{figure}
\begin{figure}
\centering \psfig{file=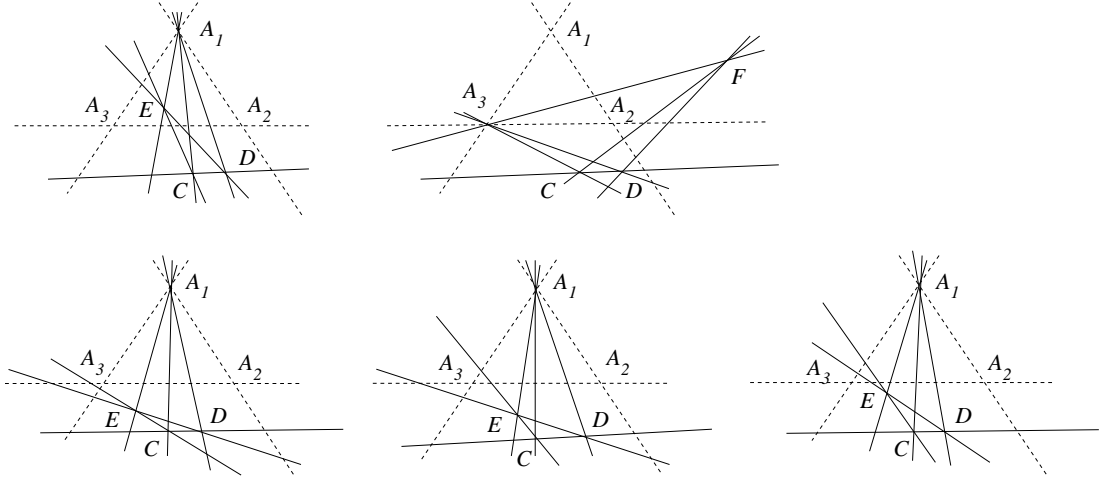}
\caption{$E \in T_0 \cup T_1$, $F \in T_3$}
\end{figure}
\begin{figure}
\centering \psfig{file=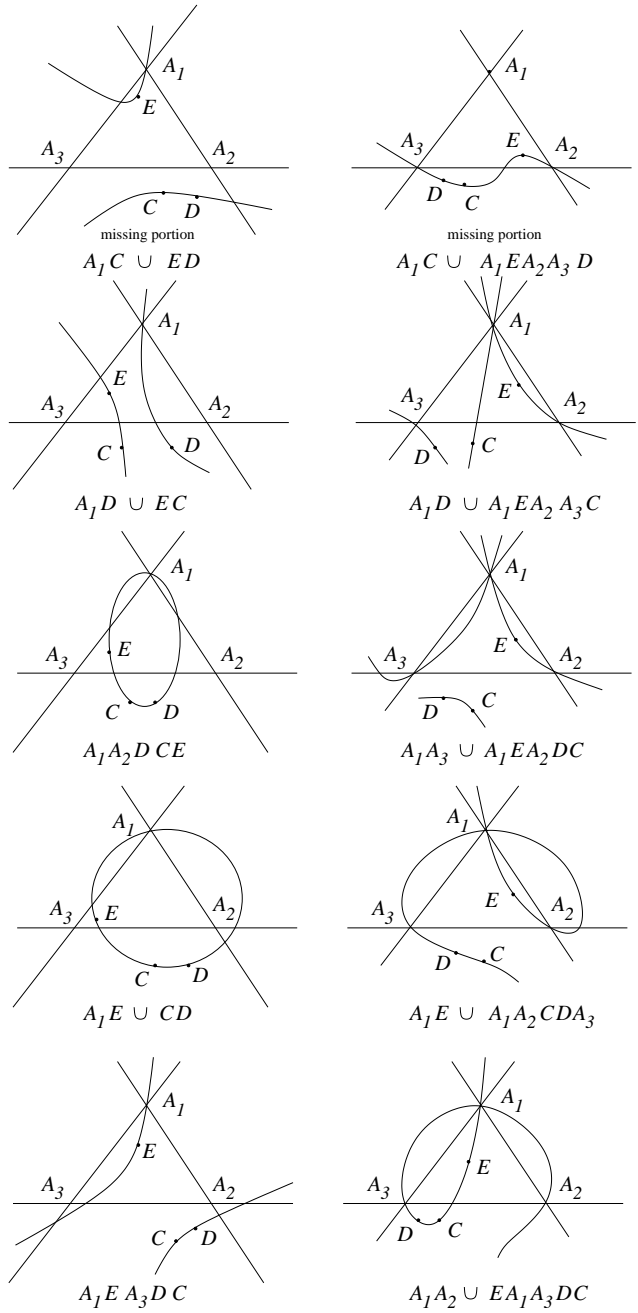}
\caption{$E \in T_0$}   
\end{figure}
\begin{figure}
\centering \psfig{file=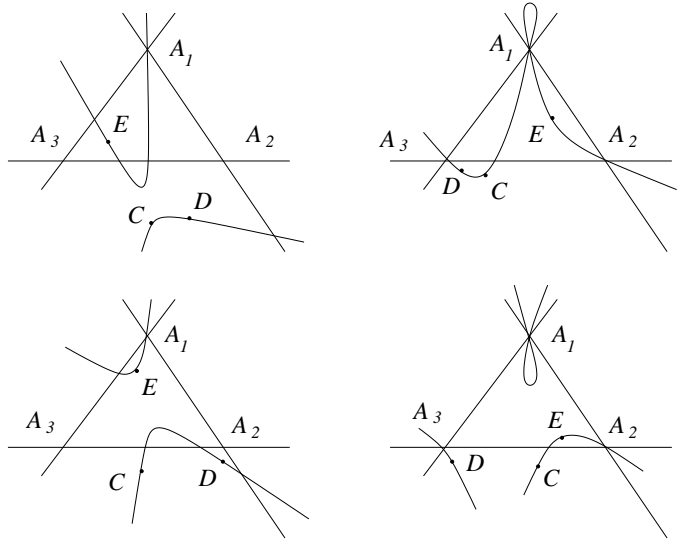}
\caption{$E \in T_0$, missing portion}   
\end{figure}
\begin{figure}
\centering \psfig{file=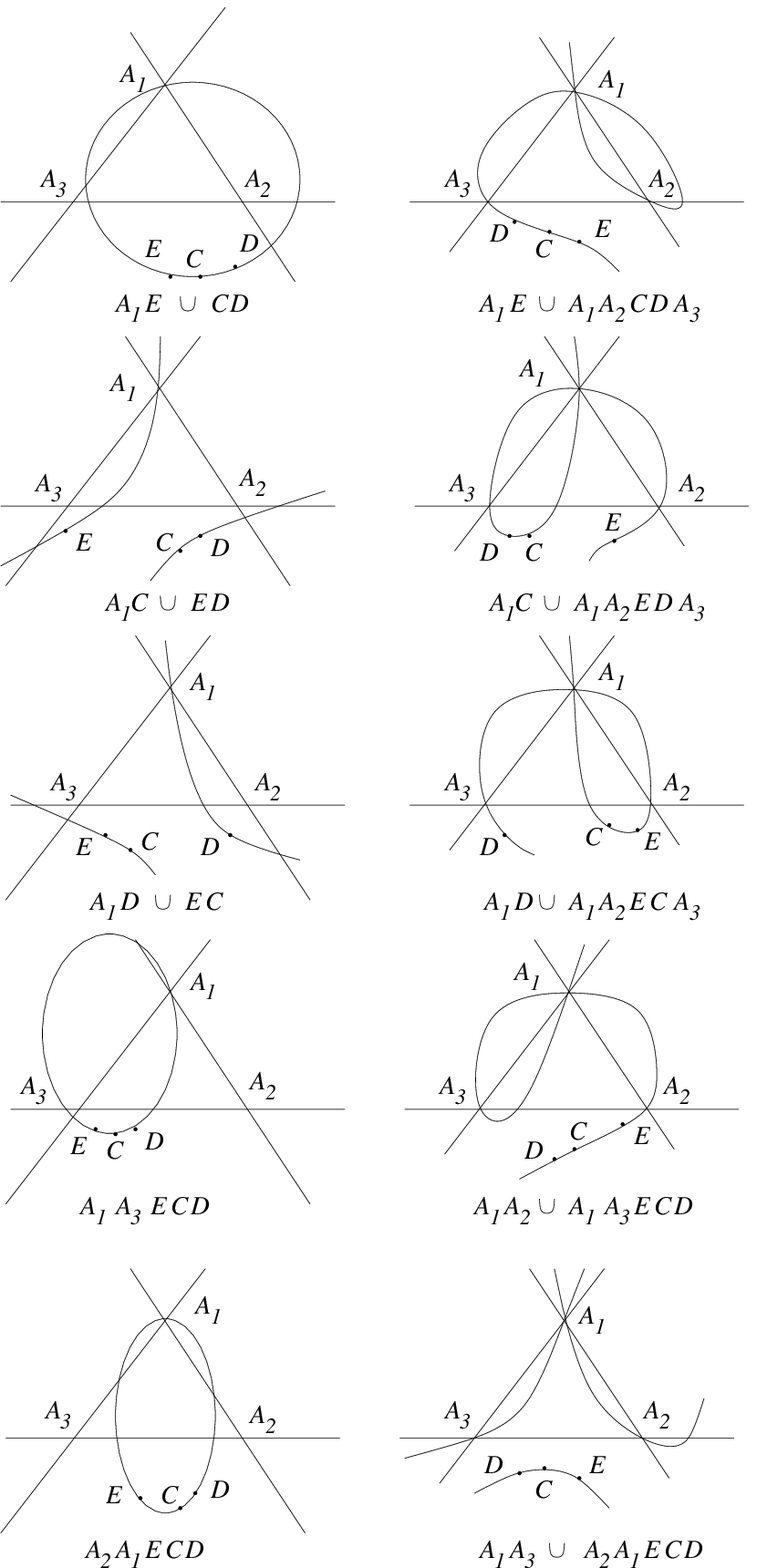}
\caption{$E \in T_1$ case 1}
\end{figure}
\begin{figure}
\centering \psfig{file=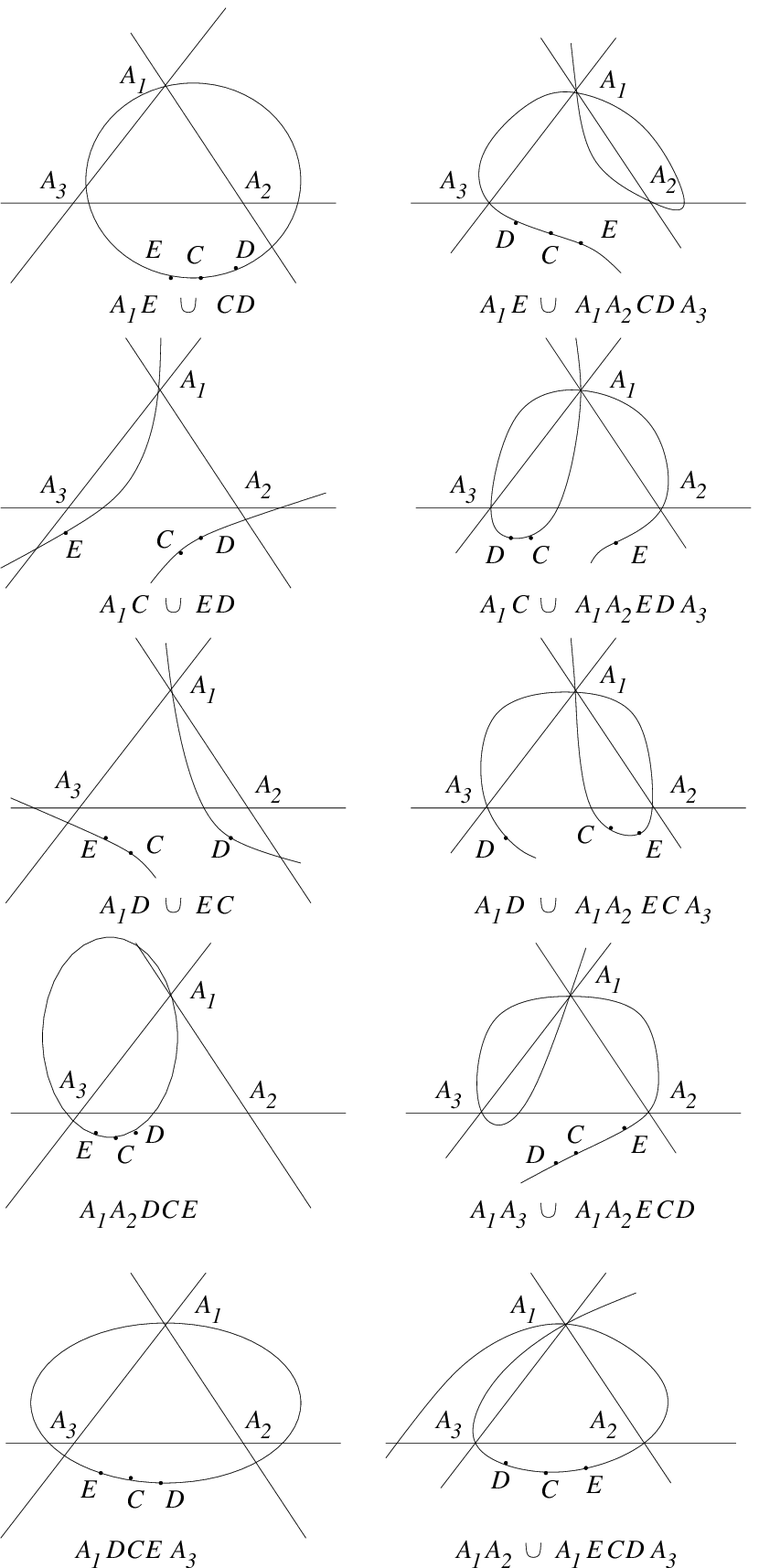}
\caption{$E \in T_1$, case 2}
\end{figure}
\begin{figure}
\centering \psfig{file=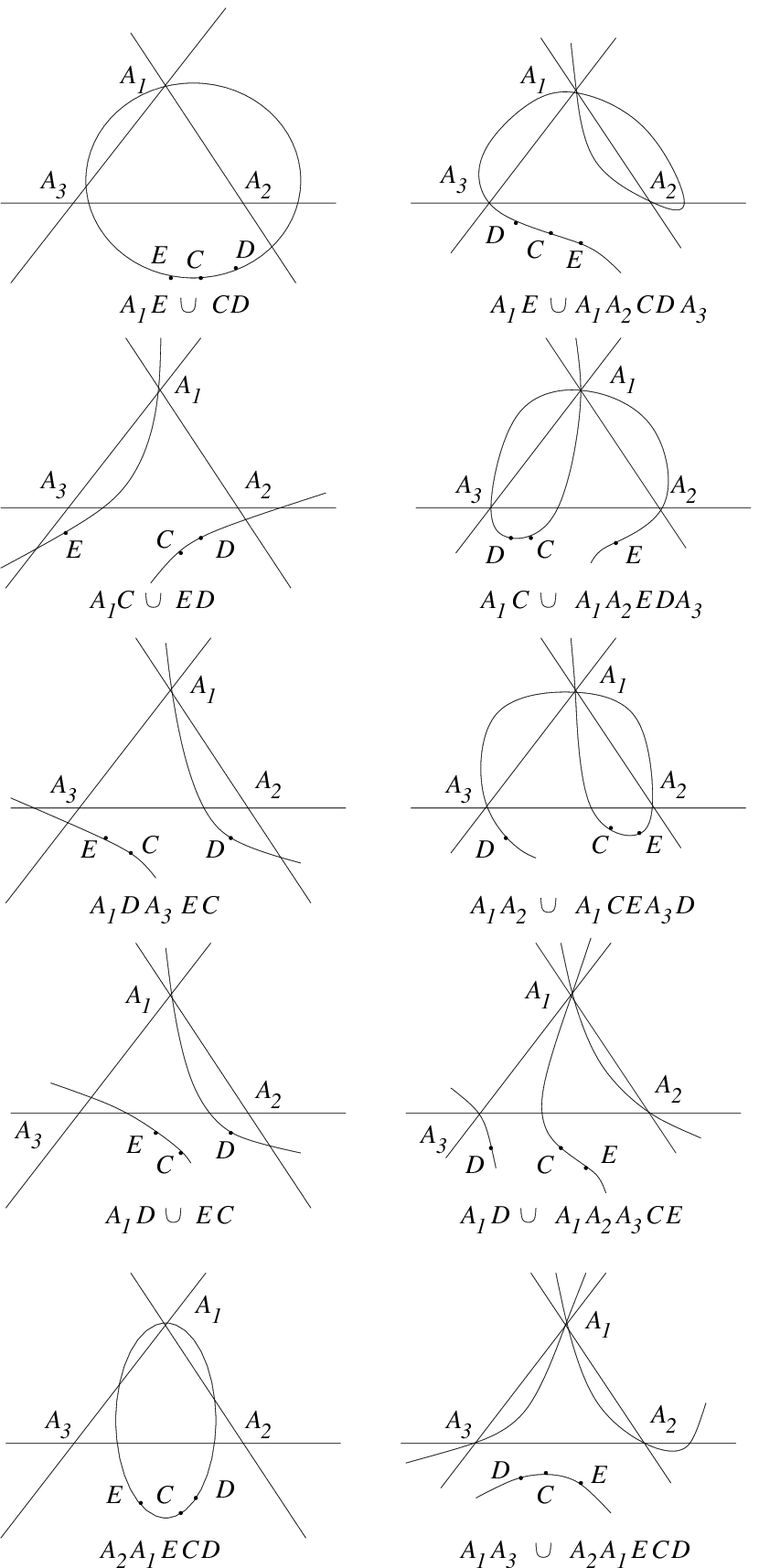}
\caption{$E \in T_1$, case 3}  
\end{figure}
\begin{figure}
\centering \psfig{file=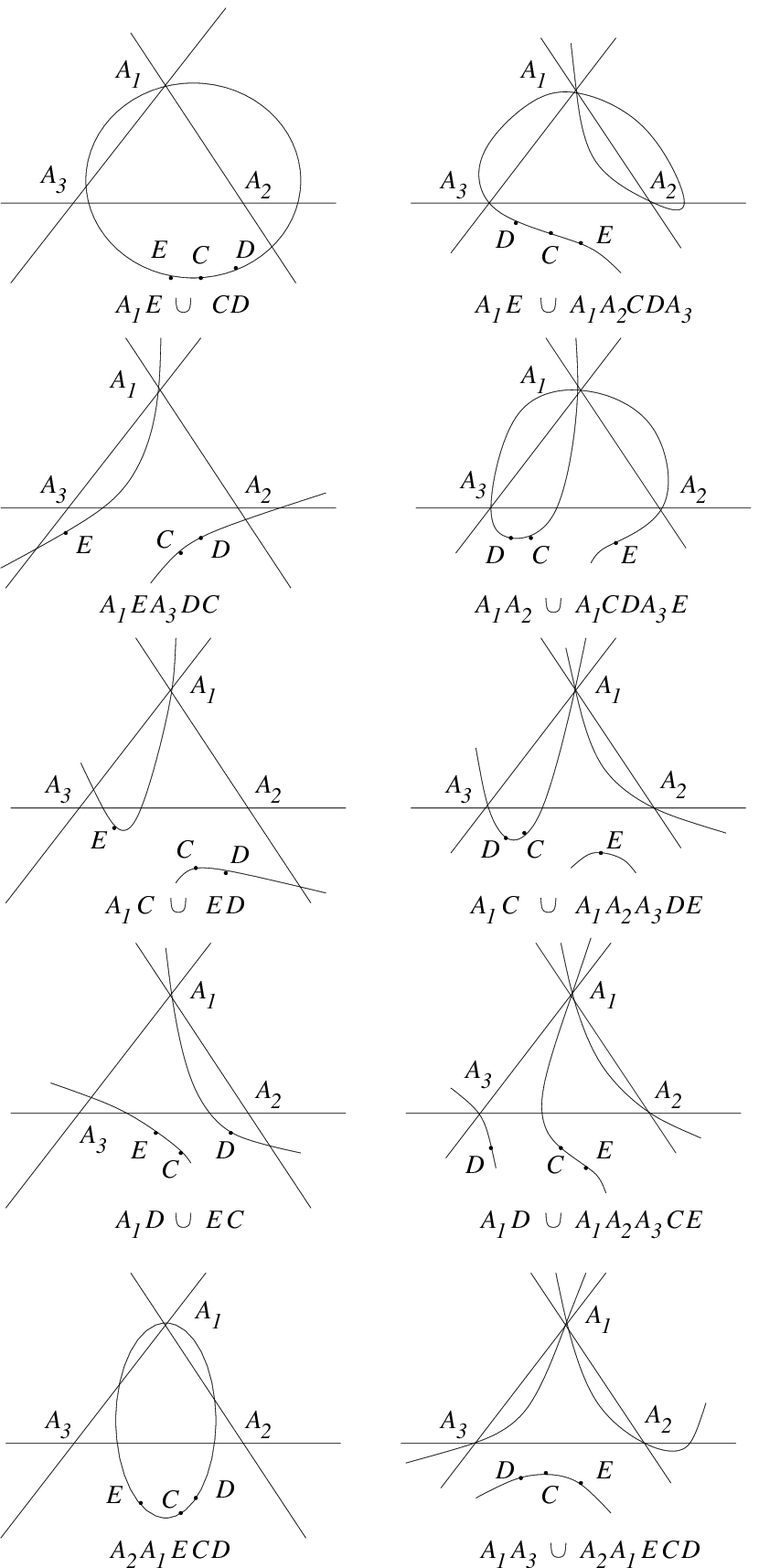}
\caption{$E \in T_1$, case 4}   
\end{figure}
\begin{figure}
\centering \psfig{file=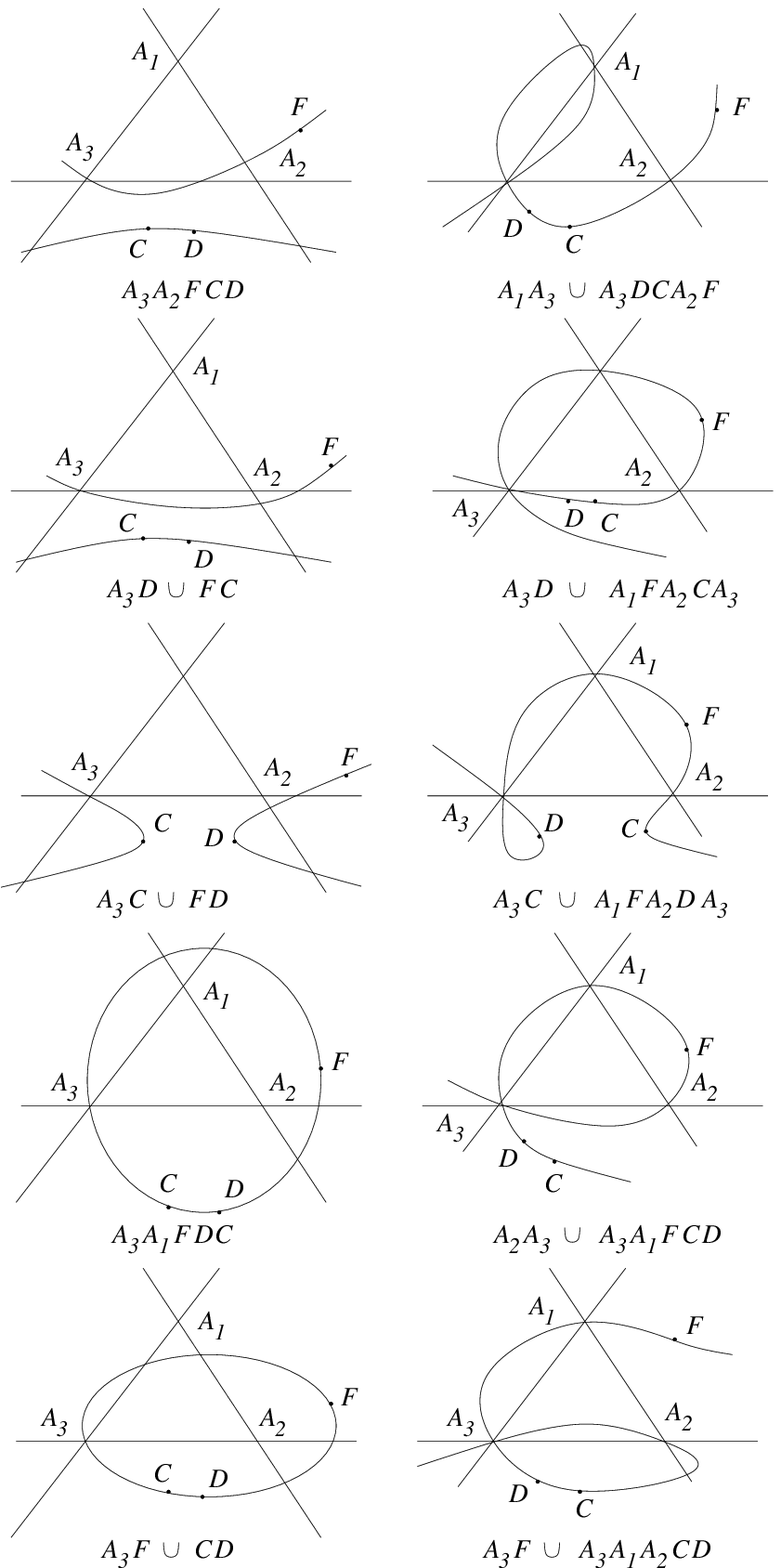}
\caption{$F \in T_3$}
\end{figure}

{\em Proof\/}
Let $O_3$ be crossing, let $E \in T_0 \cup T_1 \cup T_2$ and 
$H$ be an oval met after $E$ by the pencil of conics 
$\mathcal{F}_{A_1A_2A_3C}: A_1A_3 \cup A_2C \to A_1A_2 \cup A_3C \to A_1C \cup A_2A_3$
We shall prove that $H$ must be also in $T_0 \cup T_1 \cup T_2$. 
Assume first that $E \in T_0 \cup T_1$.

Perform the Cremona transformation $cr(A_1, A_2, A_3)$ and consider then the pencil of conics $\mathcal{F}_{A_1ECD}$. 
The possible positions for the double lines of this pencil are shown in 
Figure 8. The preimage of $\mathcal{F}_{A_1ECD}$ is the pencil of rational cubics $\mathcal{F}_{A_1A_1A_2A_3ECD}$.
This pencil has five reducible cubics, corresponding to the following conics of $\mathcal{F}_{A_1ECD}$: the three double lines and the two 
conics passing respectively through $A_2$ and $A_3$. 
The Figures 9-15 ahead show both pencils, but the actual reasoning in the proof is done with the pencil of conics.
For $E \in T_0$, there exists a conic $C_2$ in the piece  $A_1EA_3DC \to A_1C \cup ED$
of $\mathcal{F}_{A_1ECD}$, that is tangent to the base line $A_2A_3$.
The image of this conic is a cuspidal cubic.
The pencil of conics (and its image by $cr$) is shown in Figure 9. In this Figure, the portion
$C_2 \to A_1C \cup ED$ is not represented. This missing portion depends namely upon the position of the intersection of lines 
$ED \cap A_1C$: in $T_0$ or in $T_1$. Both cases are displayed in Figure 10.
If $E \in T_1$, there are four possible sequences of singular cubics (all reducible) for $\mathcal{F}_{A_1A_1A_2A_3ECD}$, 
see Figures 11-14.
Let $H$ be one of the remaining ovals of $C_{18}$. If $E \in T_0$, then $H$ must be swept out in the 
portion $\mathcal{F}_{A_1ECD}: A_1A_2DCE \to A_1E \cup CD \to A_1EA_3DC \to A_1C \cup ED$. 
Moreover, if $H$ is met after $E$ by the pencil of lines $\mathcal{F}_C: A_2 \to A_3 \to A_1$, then $H$ lies inside of 
$\mathcal{O} = cr(\mathcal{J})$.
If $E \in T_1$, then $H$ is swept out by $\mathcal{F}_{A_1ECD}$ in the portion:  

$A_2A_1ECD \to A_1E \cup CD \to A_1C \cup ED$ (case 1)

$A_1DCEA_3 \to A_1E \cup CD \to A_1C \cup ED$ (case 2)

$A_2A_1ECD \to A_1E \cup CD \to A_1C \cup ED$ (case 3)

$A_2A_1ECD \to A_1E \cup CD \to A_1EA_3DC \to A_1C \cup ED$ (case 4)

In all cases, if  $H$ is met after $E$ by the pencil of lines 
$\mathcal{F}_C: A_2 \to A_3 \to A_1$, then $H$ must lie inside of $\mathcal{O} = cr(\mathcal{J})$.   
Let now $E \in T_2$. 
Note that there are four possibilities for the pencil of rational cubics $\mathcal{F}_{A_2A_2A_1CDA_3}$, which are deduced 
from the pencils $\mathcal{F}_{A_1A_1A_2A_3ECD}, E \in T_1$ by an axial symmetry switching $(A_1, A_3)$ with $(A_2, C)$.
This finishes the case where $O_3$ is crossing. 

Let $O_3$ be non-crossing and $F \in T_3$.
Let $H$ be an oval met before $F$ by the pencil of conics $\mathcal{F}_{A_1A_2A_3C}: A_1A_3 \cup A_2C \to A_1A_2 \cup A_3C$.
We shall prove that $H$ must be also in $T_3$.
Perform the Cremona transformation $cr$, and consider
the pencil of conics $\mathcal{F}_{A_3FCD}$. The double lines of this pencil are shown in
Figure 8. The preimage of $\mathcal{F}_{A_3FCD}$ is the pencil of rational cubics $\mathcal{F}_{A_1A_2A_3A_3FCD}$.
This pencil has five singular cubics (all reducible), corresponding to the following conics of $\mathcal{F}_{A_3FCD}$: 
the three double lines and the two conics passing respectively through $A_1$ and $A_2$. 
There is only one possible sequence of singular cubics for $\mathcal{F}_{A_1A_1A_2A_3ECD}$, see Figure 15.
If $H$ is one of the remaining ovals of $C_{18}$, then $H$ must be swept out by $\mathcal{F}_{A_3FCD}$ in the portion 
$A_3A_1FDC \to A_3F \cup CD \to A_3A_2FCD$.
Moreover, if $H$ is met before $F$ by the pencil of lines $\mathcal{F}_C: A_2 \to A_3$, then $H$ must lie inside of
$\mathcal{O}$. This finishes the case where $O_3$ is non-crossing. $\Box$

\begin{lemma}
Let $C_9$ be an $M$-curve of degree $9$ with three nests and a jump.
One of the three possibilities hereafter arises:
\begin{enumerate}
\item
$\lambda_0 - \lambda_4 - \lambda_5 - \lambda_6 = 0$ and $\Pi_+-\Pi_- = 4$
\item
$O_3$ is crossing, $\lambda_0 - \lambda_4 - \lambda_5 = -1$, $\epsilon_3 = 1$, $\Pi_+-\Pi_- = 3$,
\item
$O_3$ is non-crossing, $\lambda_6 = 1$, $\epsilon_3 = -1$, $\Pi_+-\Pi_- = 3$
\end{enumerate}
\end{lemma} 

{\em Proof\/} 
Combining Lemma 17 with the fact that $\lambda_0 - \lambda_4 - \lambda_5 - \lambda_6 \leq 0$
implies that $\lambda_0 - \lambda_4 - \lambda_5 - \lambda_6 \in \{ 0, -1 \}$.
$\Box$

\section{Complex orientations again}

\subsection{Orevkov's complex orientation formulas}

Let $C_m$ be an $M$-curve of degree $m = 4d+1$, $d \geq 2$.
In this subsection, we shall call {\em nest $\mathcal{N}$ of depth n\/} a configuration of ovals 
$(o_1, \dots, o_n)$ such that $o_i$ lies in the interior of $o_j$ for all pairs $i, j$ with $j > i$.
A nest is {\em maximal\/} if it is not a subset of a bigger nest of $C_m$.
We assume that there exist four maximal nests $\mathcal{N}_i, i \in \{1, 2, 3, 4\}$ of $C_m$ verifying: if
$\mathcal{F}$ is a pencil of conics based in the four innermost ovals of the nests, any conic of $\mathcal{F}$ intersects the union of 
the four nests and $\mathcal{J}$ at at least $2m-2$ points. 
Let $V_i$ be the outermost oval of the nest $\mathcal{N}_i$.
We shall call {\em big ovals\/} the ovals that belong to the union of the nests $\mathcal{N}_i$, and 
{\em small ovals\/} the other ovals. 
For $S, s \in \{+, -\}$, let $\pi_s^S(\mathcal{N}_i)$ be the numbers of pairs of ovals $(O, o)$ with 
signs $(S, s)$ such that $O$ is an oval of $\mathcal{N}_i$ and $o$ is an empty oval contained in 
$Int(O)$. 
Let $\pi_i = (\pi_-^+ - \pi_+^+)(\mathcal{N}_i)$, $\pi'_i = (\pi_+^- - \pi_-^-)(\mathcal{N}_i)$.
Let $\Pi ^S_s(\mathcal{N}_i)$ be the number of pairs $(O, o)$ with signs $(S, s)$, where $o$ is a
small oval in $Int(O)$.
Denote by $P_i$ the number of positive ovals in $\mathcal{N}_i$, and by $N_i$ the number of non-empty
ovals among them. Denote by $Q_i$ the number of negative ovals in $\mathcal{N}_i$, and by $M_i$ the number of non-empty
ovals among them. 
Let $p_1, \dots, p_4$ be four points distributed in the innermost ovals of the four nests.
If $\mathcal{N}_i, \mathcal{N}_j, \mathcal{N}_k$ have all depth $d$, call {\em principal triangle\/} $p_ip_jp_k$
the triangle $p_ip_jp_k$ that does not intersect $\mathcal{J}$.
The formulas hereafter are proven in \cite{or6} (with slightly different notations):
  

\begin{description}
\item[First complex orientation formula (Orevkov)]
Let $C_m$ be such that the nests $\mathcal{N}_i, i \in \{1, 2, 3, 4\}$ have respective depths 
$d, d, d, d-1$, and $p_l, l \in \{1, 2, 3, 4\}$ lies in the principal triangle determined by the other three points 
$p_i, p_j, p_k$, then:   
\begin{displaymath}
\pi_i + \pi_j + \pi_k + \pi'_l = N_i^2 +  N_j^2 + N_k^2 +  M_l^2
\end{displaymath}
\end{description}

\begin{description}
\item[Second complex orientation formula (Orevkov)]
Let $C_m$ be such that: the nests $\mathcal{N}_l, l \in \{1, 2, 3, 4\}$ have all depth $d$, some $V_i, i \in \{1, 2, 3\}$ coincides 
with $V_4$, but the nests $\mathcal{N'}_i, \mathcal{N'}_4, \mathcal{N}_j, \mathcal{N}_k$ are pairwise disjoint, with $\mathcal{N}'_l = 
\mathcal{N}_l \setminus V_l$. 
Let us denote by $V$ the oval $V_i = V_4$, and by $T$ the principal triangle $p_1p_2p_3$. 
Let $\Pi_l = \Pi ^+_-(\mathcal{N}_l) - \Pi ^+_+(\mathcal{N}_l)$, 
$\Pi'_l = \Pi ^-_+(\mathcal{N}_l) - \Pi ^-_-(\mathcal{N}_l)$.
Let $Int ^+(V) = Int(V) \setminus T$, $Int ^-(V) = Int(V) \cap T$.
For any big oval $O \not= V \subset \mathcal{N}_l$, let  $Int ^\pm(O) = Int(O)$.
For $l \in \{i, 4\}$, denote by $\tilde \Pi ^S_s(\mathcal{N}_l)$ the numbers of pairs $(O, o)$ with signs $(S, s)$ where
$O \subset \mathcal{N}_l$ is big and $o \subset Int ^S(O)$ is small.
Let
$\tilde \Pi_l = \tilde \Pi ^+_-(\mathcal{N}_l) - \tilde \Pi ^+_+(\mathcal{N}_l)$,
$\tilde \Pi'_l = \tilde \Pi ^-_+(\mathcal{N}_l) - \tilde \Pi ^-_-(\mathcal{N}_l)$.
If $p_i, i \in \{1, 2, 3\}$ lies in the principal triangle $p_jp_kp_4$, then:

$\tilde \Pi'_i + \Pi_j + \Pi_k + \tilde \Pi_4 = Q^2_i - 2Q_i + P^2_j - P_j +  P^2_k - P_k +  P^2_4 - P_4 + \nu(V)$,
where $\nu(V) = 0$ if $V$ is positive, and $1$ if $V$ is negative.
\end{description}

\subsection{Proof of Theorem 1}

Let $C_9$ be an $M$-curve of degree 9 with real scheme
$\langle \mathcal{J} \amalg 1 \langle \alpha_1 \rangle
\amalg 1 \langle \alpha_2 \rangle \amalg 1 \langle \alpha_3 \rangle \amalg \beta \rangle$.
The complex scheme of $C_9$ is determined by the complex schemes of the three nests $\mathcal{O}_i = 1 \langle \alpha_i \rangle, i
\in \{1, 2, 3\}$.
The complex scheme $\mathcal{S}_i$ of a nest $\mathcal{O}_i$ is encoded as follows:
$1_{\nu_i} \langle a ^+_i \amalg a ^-_i \rangle$, where $a ^+_i +  a ^-_i = \alpha_i$, 
$a ^+_i -  a ^-_i \in \{0, \pm 1, \pm 2\}$. Let $\mu_i \in \pm$ be the sign of $a ^+_i -  a ^-_i$.
We replace the standard encoding by a simpler one, writing:
$\mathcal{S}_i = \nu_i$ if $a ^+_i -  a ^-_i = 0$,
$\mathcal{S}_i = (\nu_i, \mu_i, \mu_i)$ if $a ^+_i -  a ^-_i = \pm 2$
$\mathcal{S}_i =(\nu_i, \mu_i)$ if  $a ^+_i -  a ^-_i = \pm 1$.
Assume $O_i$ is separating and $\alpha_i$ is even.
Let $A_i$ and $A'_i$ be the two extreme interior ovals, such that $A_i \in T'_i$ and $A'_i \in T_0$.
If $A'_i$ is negative, we say that $\mathcal{O}_i$ is $(\nu_i, u)$, otherwise $\mathcal{O}_i$ is
$(\nu_i, d)$, where the letters $u$ and $d$ stand respectively for {\em up\/} and {\em down\/}.
We call {\em complex type\/} $\bar \mathcal{S}_i$ of a nest $\mathcal{O}_i$ the complex scheme of $\mathcal{O}_i$,
improved with the information whether $O_i$ is separating or not; if $O_i$ is
separating with $\alpha_i$ even, we specify whether the nest is up or down. 
If $\mathcal{S}_i = \nu_i$, then there are three possibilities for $\bar \mathcal{S}_i$: $(\nu_i, u)$,  $(\nu_i, d)$ and
$(\nu_i, n)$ ($n$ stands for {\em non-separating\/}).
If  $\mathcal{S}_i = (\nu_i, \mu_i)$, the complex type will be denoted by $(\nu_i, \mu_i, n)$ or  $(\nu_i, \mu_i, s)$
In the other cases, the oval $O_i$ is non-separating, and we use the same notation as for the
complex scheme. Let us call {\em complex type\/} of $C_9$, and 
denote by $\bar \mathcal{S}$ the triple $(\bar \mathcal{S}_1, \bar \mathcal{S}_2, \bar \mathcal{S}_3)$.

Choose base ovals $A_i, i=1, 2, 3$ and let $\mathcal{N}_i = (A_i, O_i), i=1, 2, 3$.
If $O_i$ is positive:
$\pi_i =  a_i^- - a_i^+$, $\pi'_i = 0$, $N_i = 1$, $M_i = 0$.
If $O_i$ is negative:
$\pi_i = 0$, $\pi'_i =  a_i^+ - a_i^-$, $N_i = 0$, $M_i = 1$.

\begin{lemma} 
Let $C_9$ have some exterior oval $B \in T_i, i \in \{0, 1, 2, 3\}$.
Then $E_i = 0$, where

$E_0 = \pi_1 +\pi_2 +\pi_3 - (N_1 + N_2 + N_3)$.

$E_1 = \pi'_1 +\pi_2 +\pi_3 - (M_1 + N_2 + N_3)$.

$E_2 = \pi_1 +\pi'_2 +\pi_3 - (N_1 + M_2 + N_3)$.

$E_3 = \pi_1 +\pi_2 +\pi'_3 - (N_1 + N_2 + M_3)$.
\end{lemma}

{\em Proof\/}  
The first formula applies making $\mathcal{N}_i, i=1, 2, 3$ and $\mathcal{N}_4 = B$. $\Box$
 
Let $O_i$ be separating. 
Again, choose base ovals $A_1, A_2, A_3$, and let $A_4$ be a fourth oval, interior to $O_i$, lying in 
$T_i$. Let  $\mathcal{N}_l = (A_l, O_l), l=1, 2, 3$ and $\mathcal{N}_4 = (A_4, O_i)$.
Let: $F_i = \tilde \Pi'_i + \tilde \Pi_4 - (Q^2_i - 2Q_i + P^2_4 - P_4 + \nu(O_i))$, and
$G_l =  P^2_l - P_l -  \Pi_l$. It is easily seen that $G_l$ depends only on $\mathcal{S}_i$, and
$F_i$ depends only on $\bar \mathcal{S}_i$.

\begin{lemma}
Let $C_9$ have three nests and separating $O_i$. Then, 
$F_i = G_j + G_k$
\end{lemma}

{\em Proof\/} The second formula applies with the nests $\mathcal{N}_l, l=1, 2, 3, 4$. $\Box$

In Figure 16, 17 we computed the terms appearing in Lemmas 19, 20.

\begin{lemma}
Let $C_9$ be an $M$-curve of degree 9 with three nests.
If the union $T_0 \cup T_1 \cup T_2 \cup T_3$ is empty, then 
$C_9$ verifies: $\mathcal{S}_1, \mathcal{S}_2 \in \{(+, -), (-, +)\}$, $S_3 \in \{(+, -, -), (-, +, +)\}$.
\end{lemma}

{\em Proof\/} One has $\lambda_0 - \lambda_4 - \lambda_5 - \lambda_6 = 0$, $\Lambda_+ - \Lambda_- = 0$,
$\Pi_+ - \Pi_- = 4$. $\Box$
 
\begin{theorem}
Let $C_9$ be an $M$-curves of degree 9 with real scheme $\langle \mathcal{J} \amalg 1 \langle \alpha_1 \rangle
\amalg 1 \langle \alpha_2 \rangle \amalg 1 \langle \alpha_3 \rangle \amalg \beta \rangle$. At least one of the
$\alpha_i, i=1, 2, 3$ is odd.
\end{theorem}

{\em Proof\/} 
Assume that $\alpha_1, \alpha_2, \alpha_3$ are even.
By Lemma 18, $C_9$ has no jump, so this curve realizes one of the four complex schemes of Figure 18.   
The last column $\mathcal{Z}$ of this figure contains the indices of the triangles $T_i, i \in \{0, 1, 2, 3\}$
that may contain exterior ovals. 
In Figure 19, we assume that $O_i$ is separating and compute the term $F_i - G_j - G_k$.
In Figure 20, we list the a priori possible complex types for $C_9$ along with the data $\mathcal{Z}$.
The first two types contradict Lemma 21. 
For the other complex types, choose each $A_i, i \in \{1, 2, 3\}$ in such a way that $(A_i, O_i)$ is a positive pair. 
The last identity in Lemma 10 reads
$\lambda_0 - \lambda_4 - \lambda_5 - \lambda_6 = -4$.
Compute the values of $\lambda_0, \lambda_4, \lambda_5, \lambda_6$.
For the last four types, one gets $\lambda_0 \leq -4$, which contradicts Lemma 16.
For the remaining two types, one has $\lambda_6 = 4$ or $5$, which contradicts Proposition 1. $\Box$ 

\begin{figure}
\centering
\begin{tabular}{c c c c c c}
$\mathcal{S}_l$ & $\pi_l$ & $\pi'_l$ & $N_l$ & $M_l$ & $G_l$\\
$-$ & 0 & 0 & 0 & 1  & 0\\ 
$+$ & 0 & 0 & 1 & 0  & 1\\
$(-, +)$ & 0 & 1 & 0 & 1 & 0\\
$(+, -)$ & 1 & 0 & 1 & 0 & 0\\
$(-, -)$ & 0 & $-1$ & 0 & 1 & 0\\
$(+, +)$ & $-1$ & 0 & 1 & 0 & 2\\
$(-, +, +)$ & 0 & 2 & 0 & 1 & 0\\
$(+, -, -)$ & 2 & 0 & 1 & 0 & $-1$\\
$(-, -, -)$ & 0 & $-2$ & 0 & 1 & 0\\
$(+, +, +)$ & $-2$ & 0 & 1 & 0 & 3
\end{tabular}
\caption{}
\end{figure}

\begin{figure}
\centering
\begin{tabular}{c c c c c}
$\bar \mathcal{S}_i$ & $F_i$ & & $\bar \mathcal{S}_i$ & $F_i$\\
$(-, d)$ & 0 & & $(-, -, s)$ & -1\\
$(-, u)$ & -1 & & $(-, +, s)$ & 0\\
$(+, d)$ & 0 & & $(+, -, s)$ & 0\\
$(+, u)$ & -1 & & $(+, +, s)$ & -1\\
\end{tabular}
\caption{}
\end{figure}

\begin{figure}
\centering
\begin{tabular}{c c c c c c c c}
$\mathcal{S}_1$ & $\mathcal{S}_2$ & $\mathcal{S}_3$ & $E_0$ & $E_1$ & $E_2$ & $E_3$ & $\mathcal{Z}$\\
$-$ & $-$ & $-$ & $0$ & $-1$ & $-1$ & $-1$ & (0)\\
$-$ & $-$ & + & $-1$ & $-2$ & $-2$ & $0$ & (3)\\
$-$ & + & + & $-2$ & $-3$ & $-1$ & $-1$ & $\emptyset$\\
+ & + & + & $-3$ & $-2$ & $-2$ & $-2$ & $\emptyset$
\end{tabular}
\caption{}
\end{figure}

\begin{figure}
\centering
\begin{tabular}{c c c c}
$\bar \mathcal{S}_i$ & $\mathcal{S}_j$ & $\mathcal{S}_k$ & $F_i - G_j - G_k$\\
$(-, d)$ & $-$ & $-$ & $0$\\
$(-, d)$ & $-$ & $+$ & $-1$\\
$(-, d)$ & $+$ & $+$ & $-2$\\
$(-, u)$ & $-$ & $-$ & $-1$\\
$(-, u)$ & $-$ & $+$ & $-2$\\
$(-, u)$ & $+$ & $+$ & $-3$\\
$(+, d)$ & $-$ & $-$ & $0$\\
$(+, d)$ & $-$ & $+$ & $-1$\\
$(+, d)$ & $+$ & $+$ & $-2$\\
$(+, u)$ & $-$ & $-$ & $-1$\\
$(+, u)$ & $-$ & $+$ & $-2$\\
$(+, u)$ & $+$ & $+$ & $-3$
\end{tabular}
\caption{}
\end{figure}

\begin{figure}
\centering
\begin{tabular}{c c c c}
$\bar \mathcal{S}_1$ & $\bar \mathcal{S}_2$ & $\bar \mathcal{S}_3$ & $\mathcal{Z}$\\
$(+, n)$ & $(+, n)$ & $(+, n)$ & $\emptyset$\\
$(-, n)$ & $(+, n)$ & $(+, n)$ & $\emptyset$\\
$(-, n)$ & $(-, n)$ & $(+, n)$ & $(3)$\\
$(-, n)$ & $(-, n)$ & $(+, d)$ & $(3)$\\
$(-, n)$ & $(-, n)$ & $(-, n)$ & $(0)$\\ 
$(-, d)$ & $(-, n)$ & $(-, n)$ & $(0)$\\ 
$(-, d)$ & $(-, d)$ & $(-, n)$ & $(0)$\\ 
$(-, d)$ & $(-, d)$ & $(-, d)$ & $(0)$
\end{tabular}
\caption{}
\end{figure} 

\subsection{More on complex orientations and rigid isotopy}

\begin{proposition}
Let $C_9$ be an $M$-curve with three nests, and $A_1, A_2, A_3$ be base ovals such that $T_0$
contains only exterior ovals. Then $\vert \lambda_0 \vert \leq 2$.
\end{proposition}

{\em Proof\/} Let $C_9$ satisfy the conditions of the proposition and assume that $\vert \lambda_0 \vert > 2$.
By Lemma 16, one has $\lambda_0 = \pm 3$, $\epsilon_i = \mp 1$ for $i = 1, \dots, 6$, and the non-empty ovals are all
three separating.
At least one of the nests contains an odd numer of ovals by Theorem 1. 
The a priori possible complex schemes are listed in Figure 21, along with the value of $E_0$.
By Lemma 19,  one must have $E_0 = 0$, this excludes three complex schemes.
In Figure 22, we compute the values of $E_1, E_2, E_3$ for the remaining cases.
By Lemma 19, the triangles $T_1, T_2, T_3$ contain no exterior ovals.
By Lemma 16, one of the quadrangles, say $Q_1$, is empty.
Applying the first identity of Lemma 10, with $\lambda_0 = 3$, $\lambda_1 = 0$  
yields $\lambda_4 = 1$. Hence, the complex scheme $S_1$ is $-$.
The complex type $\bar S_1$ is $(-, u)$. For both cases $S_2, S_3 = -, (-, -)$ or $(-, -), (-, -)$, one has
$F_1 - G_2 - G_3 = -1$. This contradicts Lemma 20. $\Box$

\begin{figure}
\centering
\begin{tabular}{ c c c c}
$S_1$ & $S_2$ & $S_3$ & $E_0$\\
$-$ & $-$ & $(-, -)$ & $0$\\
$-$ & $(-, -)$ & $(-, -)$ & $0$\\
$(-, -)$ & $(-, -)$ & $(-, -)$ & $0$\\
$+$ & $+$ & $(+, +)$ & $-2$\\
$+$ & $(+, +)$ & $(+, +)$ & $-5$\\
$(+, +)$ & $(+, +)$ & $(+, +)$ & $-6$\\
\end{tabular}
\caption{}
\end{figure}  

\begin{figure}
\centering
\begin{tabular}{ c c c c c c }
$S_1$ & $S_2$ & $S_3$ & $E_1$ & $E_2$ & $E_3$\\
$-$ & $-$ & $(-, -)$ & $-1$ & $-1$ & $-2$\\
$-$ & $(-, -)$ & $(-, -)$ & $-1$ & $-2$ & $-2$\\
$(-, -)$ & $(-, -)$ & $(-, -)$ & $-2$ & $-2$ & $-2$\\
\end{tabular}
\caption{}
\end{figure}

\begin{proposition}
Let $C_9$ be an $M$-curve with three nests and a jump arising from $B, C, D$ where $B, C$ are the extreme ovals in
the chain of $O_3$, and $A_1, A_2, C, D, B$ lie in convex position. 
Assume there exist front and back ovals. Then the chain 
of $O_3$ splits into two consecutive subchains
$(G_i), (H_j)$, where the $G_i$ are front and the $H_j$ are back.
\end{proposition}  

{\em Proof\/} 
Consider two supplementary ovals $E, F$ such that $B, D, E, F, C$ are disposed
in this ordering in the chain of $O_3$.
Assume that
\begin{enumerate}
\item
$D, F$ are front and $E$ is back, or
\item
$E$ is front and $D, F$ are back.
\end{enumerate}

Perform the Cremona transformation $cr$ and consider the pencil of conics $\mathcal{F}_{A_2CDE}$.
The double lines of this pencil for cases 1 and 2 are shown in Figure 23, 
along with the possible positions of $F$.
In both cases, $F$ is swept out in the portions $CE \cup A_2D \to CD \cup A_2E \to CA_2 \cup DE$.
Note that the base points $A_1$ and $A_3$ are swept out outside of this portion.
Perform the Cremona transformation back, the conic through $C, A_2, D, E, F$ is mapped onto a rational cubic $C_3$
passing through $A_1, A_2, A_3, C, D, E, F$, with node at $A_2$.
Figure 24 shows the pair of possible cubics $C_3$ for either case 1 and 2.
Any one of the four cubics intersects $C_9$ at $29$ points, this is a contradiction. $\Box$    

\begin{figure}
\centering \psfig{file=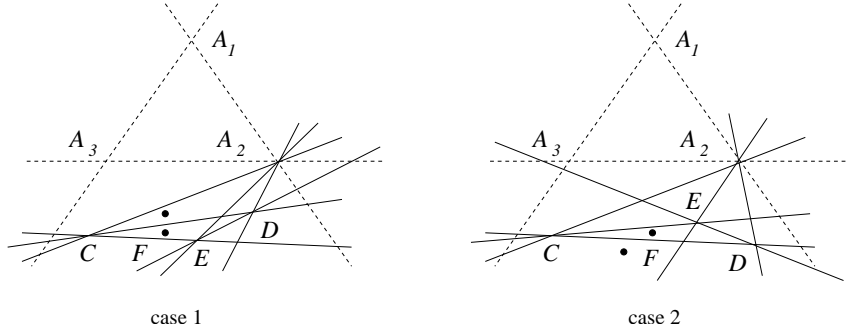}
\caption{The double lines}
\end{figure}

\begin{figure}
\centering \psfig{file=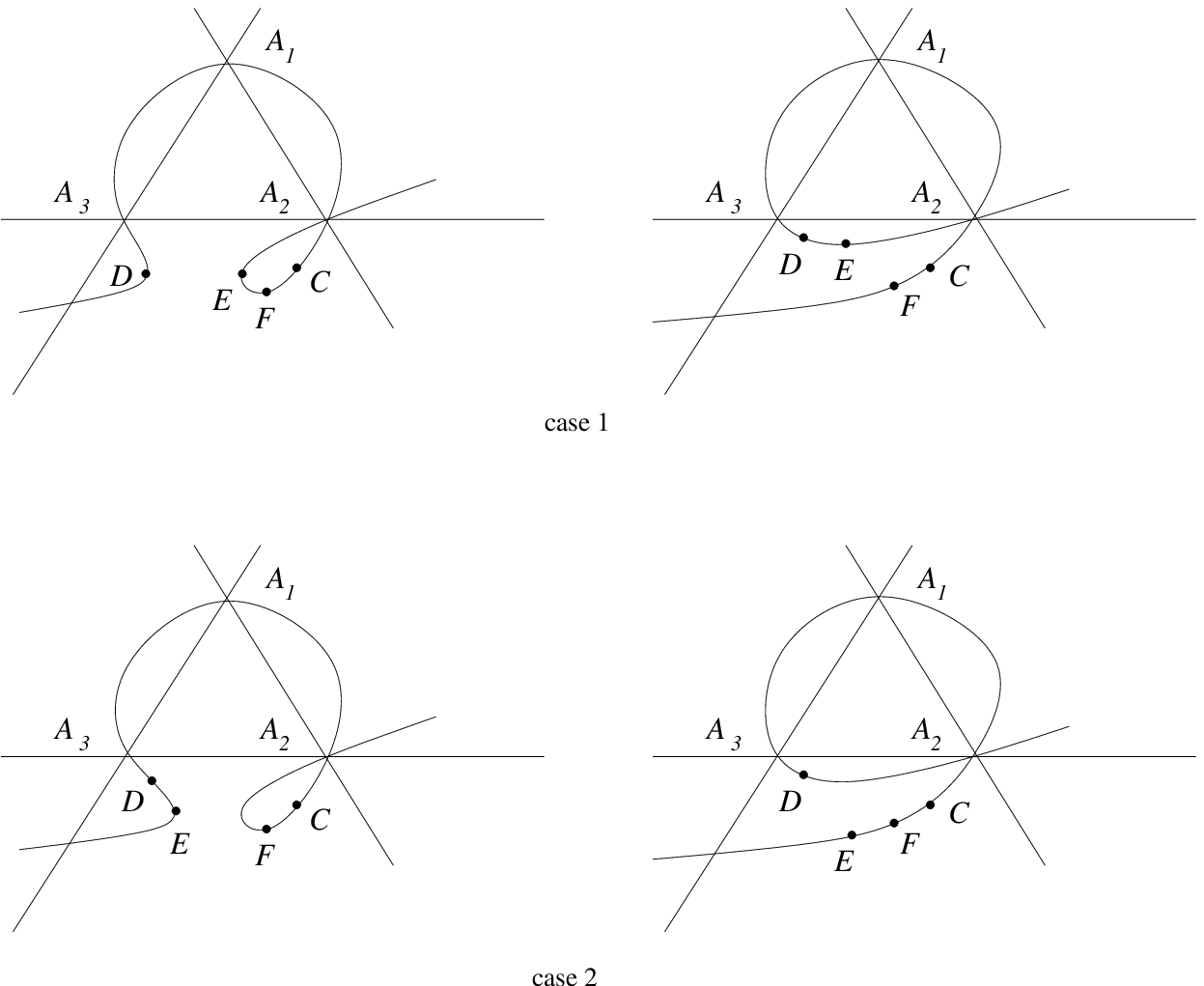}
\caption{The possible cubics $C_3$}
\end{figure}

\vspace{1cm}
{\sc Acknowledgement} I am grateful to Stepan Orevkov for his useful remarks.

\vspace{3cm}
severine.fiedler@live.fr

\end{document}